\documentclass[11pt]{article}
\usepackage{amsmath}
\usepackage{amsfonts}
\usepackage{amssymb}
\usepackage{bbm}
\usepackage{latexsym}
\usepackage[dvips]{epsfig}
\usepackage{graphicx}
\usepackage{psfrag}
\usepackage{pstricks}
\usepackage{colordvi}
\usepackage{wrapfig}
\usepackage{varioref}
\usepackage{hyperref}
\usepackage{cleveref}
\usepackage{tikz}
\usepackage{pgfplots}

\usepackage{mathrsfs}
\usepackage{upgreek}

 \let\oldsection\section
\renewcommand{\section}{
  \renewcommand{\theequation}{\thesection.\arabic{equation}}
  \oldsection}

\usepackage{caption}
\usepackage{subcaption}
\usepackage[margin=0.25cm]{caption}
\usepackage{cancel}
\usepackage{pdfpages}
\usepackage[font=small,skip=0pt]{caption}

\graphicspath{{Figures/}} 

\makeatletter\@addtoreset{equation}{section}\makeatother
\newcounter{constantsnumber}

\newcommand{\bcM}{\mbox{\boldmath{$\mathcal{M}$}}}

\newcommand{\ba}{\mbox{\boldmath{$a$}}}

\newcommand{\bB}{\mbox{\boldmath{$B$}}}

\newcommand{\bd}{\mbox{\boldmath{$d$}}}
\newcommand{\bD}{\mbox{\boldmath{$D$}}}

\newcommand{\be}{\mbox{\boldmath{$e$}}}
\newcommand{\fb}{\mbox{\boldmath{$f$}}}

\newcommand{\bh}{\mbox{\boldmath{$h$}}}

\newcommand{\bg}{\mbox{\boldmath{$g$}}}
\newcommand{\bI}{\mbox{\boldmath{$I$}}}
\newcommand{\bK}{\mbox{\boldmath{$K$}}}

\newcommand{\bM}{\mbox{\boldmath{$M$}}}

\newcommand{\bn}{\mbox{\boldmath{$n$}}}
\newcommand{\bN}{\mbox{\boldmath{$N$}}}

\newcommand{\bs}{\mbox{\boldmath{$s$}}}

\newcommand{\bu}{\mbox{\boldmath{$u$}}}
\newcommand{\bU}{\mbox{\boldmath{$U$}}}

\newcommand{\bv}{\mbox{\boldmath{$v$}}}

\newcommand{\bx}{\mbox{\boldmath{$x$}}}

\newcommand{\bepsilon}{\mbox{\boldmath{$\varepsilon$}}}

\newcommand{\bvarphi}{\mbox{\boldmath{$\varphi$}}}
\newcommand{\bsigma}{\mbox{\boldmath{$\sigma$}}}

\newcommand{\BBC}{\mbox{$\mathbb{C}$}}

\newcommand{\bzero}{\mbox{$\bf 0$}}

\newcommand{\BBR}{\mbox{$\mathbb{R}$}}

\newfont{\twelvemsb}{msbm10 at 11.6pt}

\newcommand{\half}{\mbox{$\frac{1}{2}$}}

\newcommand{\tr}{\mathop{\rm tr}}

\newcommand{\qed}{\nobreak \ifvmode \relax \else
      \ifdim\lastskip<1.5em \hskip-\lastskip
      \hskip1.5em plus0em minus0.5em \fi \nobreak
      \vrule height0.75em width0.5em depth0.25em\fi}

\begin{document}
	\begin{center}
		{\bf \Large A Virtual Element Method for transversely isotropic elasticity}\vspace{2ex}\\
		D. van Huyssteen\footnote{Department of Mechanical Engineering and Centre for Research in Computational and Applied Mechanics, University of Cape Town, 7701 Rondebosch, South Africa. Email: vhydan001@uct.ac.za}
		\qquad
		B. D. Reddy\footnote{Department of Mathematics and Applied Mathematics and Centre for Research in Computational and Applied Mechanics, University of Cape Town, 7701 Rondebosch, South Africa. Email: daya.reddy@uct.ac.za}
		\\ 
	\end{center}
	
	\bigskip

	\paragraph*{\textbf{Abstract}}
	This work studies the approximation of plane problems concerning transversely isotropic elasticity, using a low-order Virtual Element Method (VEM), with a focus on near-incompressibility and near-inextensibility. Additionally, both homogeneous problems, in which the plane of isotropy is fixed; and non-homogeneous problems, in which the fibre direction defining the isotropy plane varies with position, are explored. In the latter case various options are considered for approximating the non-homogeneous fibre directions at element level. Through a range of numerical examples the VEM approximations are shown to be robust and locking-free for several element geometries and for fibre directions that correspond to mild and strong non-homogeneity.
	
	%
	
	\section{Introduction}
	The popular finite element method has the status of  a classical approach for obtaining approximate solutions to problems formulated as systems of partial differential equations or inequalities, or alternatively, in their variational form. Particularly in the domains of solid and fluid mechanics, the method has been used with great success for problems with high degrees of complexity such as non-linear problems and problems with intricate geometries (see for example the treatments in \cite{Belytschko-Liu-Moran2000,Wriggers2008}). 
	
	A number of variants of standard conforming finite element methods have been developed over the last four decades, with a range of motivations in mind. Mixed methods, for example, have allowed all variables of interest to be approximated explicitly; and in addition have provided avenues through which stable and convergent finite element approximations can be developed in situations where the selection of the values of certain parameters might lead to non-convergence. Key examples are those of near-incompressibility, or problems in structural mechanics in which the geometry is characterized by a small length scale. These two features lead, in the context of low-order standard finite element methods, to volumetric and shear locking respectively. Phenomena that may be circumvented by the use of mixed methods \cite{Boffi-Brezzi-Fortin2013,Hughes1987}.
	
	Yet another variant of the standard conforming finite element method is the discontinuous Galerkin (DG) method, in which interelement continuity is abandoned in favour of greater flexibility with regard to meshing (see for example \cite{Arnold-Brezzi-Cockburn-Marini2002}). In addition, DG methods, when designed appropriately, are stable and uniformly convergent in situations of near-incompressibility for low-order approximations \cite{Grieshaber-McBride-Reddy2015,Hansbo-Larson2002,Wihler2004}.
	
	A more recent development in the context of finite element methods is the Virtual Element Method (VEM). In contrast to the geometric restrictions on finite elements, which are most generally triangular or quadrilateral in 2D, and tetrahedral or hexahedral in 3D, the VEM permits elements to be  arbitrary polygons in 2D or polyhedra 3D. Furthermore, there is no need for elements to be convex, and degeneracies such as element sides having small interior angles or arbitrarily small edges pose no problems. Some key works in a rapidly growing literature include 
	\cite{Beirao-etal2013a,Beirao-etal2014,Gain-etal2014}. Applications of the VEM to nonlinear problems include works on nonlinear elasticity \cite{Chi-Beirao-Paulino2017,Wriggers-Reddy-Rust-Hudobivnik2017},
	elastoplasticity \cite{Artioli-etal2017b,Beirao-etal2015,Wriggers-Hudobivnik2017}, and contact \cite{Wriggers-Rust-Reddy2016}.
	
	Applications of the VEM to elasticity have been largely confined to the isotropic problem, although there have been treatments of inextensible materials \cite{Wriggers-Hudobivnik-Korelc2018}. Problems involving anisotropic materials pose additional challenges in the context of VEM approaches, particularly for non-homogeneous materials in which the anisotropy varies with position. In \cite{Auricchio-Scalet-Wriggers2017} limiting extensibility is investigated in an otherwise isotropic material using penalty, Lagrange multiplier, and perturbed Lagrangian approaches. The work was considerably extended by Rasalofson et. al. \cite{Rasolofoson-Grieshaber-Reddy2018}. This work presents a detailed treatment of the boundary value problem for transversely isotropic linear elastic materials. Conditions for well-posedness are established, and finite element approximations are studied using both conforming and reduced integration approaches. An error analysis gives an indication of conditions under which low-order approximations are uniformly convergent in the incompressible and inextensible limits, and a set of numerical experiments provides further insight into the conditions under which locking-free behaviour occurs. Specifically, with the degree of anisotropy measured through the ratio of Young's modulus in the fibre direction relative to that in the plane of isotropy, it is shown in this work that locking-free behaviour occurs in conditions of moderate anisotropy for low-order conforming quadrilaterals, in contrast to the situation for isotropic materials. Furthermore, for high degrees of anisotropy leading to near-inextensibility, locking occurs, but is circumvented by the use of selective under-integration. 
	
	The purpose of this work is to study low-order VEM approximations for plane problems concerning transversely isotropic elasticity. Of particular interest is the behaviour of VEM approximations for the limiting situations of near-incompressibility and near-inextensibility. Whereas in the case of conventional finite element approximations, as discussed above, some form of modification such as selective under-integration is necessary in order to circumvent locking, in the case of VEM approximations locking-free behaviour is observed in the incompressible and inextensible limits. 
	
	A further novel aspect of this work is the treatment of non-homogeneous transverse isotropy; that is, situations in which fibre directions vary with position. Here it becomes necessary to approximate the non-homogeneous terms appropriately in order to preserve the simplicity of the VEM approach, in which integrals are evaluated only on element boundaries. The approximations adopted are shown to be robust, with the locking-free behaviour also evident for the non-homogeneous problem.  
	
	The structure of the rest of this work is as follows. Section 2 sets out the details of the constitutive relations for transversely isotropic linear elastic materials, the set of governing equations, and the associated weak formulation. The details of the Virtual Element Method are presented in Section 3, and the set of numerical results are presented and discussed in Section 4. This work concludes with a summary of results and a discussion of open problems. 
	
	\section{The governing equations for transverse isotropy}
	Consider a linear elastic body which occupies a plane, polygonal bounded domain $\Omega \subset \BBR^2$ with boundary $\partial\Omega$. The boundary comprises a non-trivial Dirichlet part $\Gamma_D$ and Neumann part $\Gamma_N$ such that $\Gamma_D \cap \Gamma_N = \emptyset$ and $\overline{\Gamma_D \cup \Gamma_N} = \partial \Omega$.
	
	\subsection{The elastic relation}
	\label{subsec:ElasticRelation}
	Transversely isotropic materials exhibit isotropic behaviour in a specified plane, this plane being defined by a normal vector, and referred to also as the fibre direction. 
	
	The Cauchy stress tensor $\bsigma$ is related to the infinitesimal strain tensor $\bepsilon$ through the elastic relation
	\begin{equation}
	\bsigma = \BBC\bepsilon.
	\label{stress-strain}
	\end{equation}
	Here $\BBC$ is a fourth-order tensor of elastic moduli. For a transversely isotropic material with the direction of transverse isotropy defined by the unit vector $\ba$, \eqref{stress-strain} takes the form \cite{Rasolofoson-Grieshaber-Reddy2018}
	\begin{align}
	\bsigma &= \lambda(\tr\bepsilon)\bI + 2\mu_T \bepsilon + \beta (\bM:\bepsilon)\bM + \alpha ((\bM:\bepsilon)\bI + (\tr\bepsilon)\bM) + 2(\mu_L - \mu_T) (\bepsilon\bM + \bM\bepsilon).
	\label{stress-strain2}
	\end{align}
	Here $\bM = \ba\otimes \ba$, $\lambda$ and $\mu_T$ are the conventional Lam\'e parameters, $\mu_L$ is the shear modulus in the fibre direction, and $\bI$ denotes the second-order identity tensor. The material constants $\alpha$ and $\beta$ do not have a direct interpretation, though it will be seen that $\beta \rightarrow \infty$ in the limit of inextensible behaviour in the fibre direction.
	
	The special case of an isotropic material is recovered by setting $\alpha = \beta = 0$ and $\mu_L = \mu_T$.
	
	The five material constants in \eqref{stress-strain2} may be related to the ``engineering" constants, viz. Young's moduli $E_L$ and $E_T$ in the fibre direction and plane of isotropy, respectively, and the corresponding Poisson's ratios $\nu_L$ and $\nu_T$, by inverting \eqref{stress-strain2}, specializing it to the case in which $\ba =\be_3$, and comparing with the compliance relation written in the form (see for example \cite{Exadaktylos2001})
	\begin{equation}
	\label{strain_stress_engineer_csts}
	\begin{pmatrix}
	\varepsilon_{11} \\ \varepsilon_{22} \\ \varepsilon_{33} \\ 2\epsilon_{23} \\ 2\epsilon_{13} \\ 2\epsilon_{12}
	\end{pmatrix}
	= 
	\begin{pmatrix}
	\dfrac{1}{E_T} & -\dfrac{\nu_T}{E_T} & -\dfrac{\nu_L}{E_L} & 0 & 0 & 0\\
	-\dfrac{\nu_T}{E_T} & \dfrac{1}{E_T} & -\dfrac{\nu_L}{E_L} & 0 & 0 & 0\\
	-\dfrac{\nu_L}{E_L} & -\dfrac{\nu_L}{E_L} & \dfrac{1}{E_L} & 0 & 0 & 0\\
	0 & 0 & 0 & \dfrac{1}{\mu_L} & 0 & 0\\
	0 & 0 & 0 & 0 & \dfrac{1}{\mu_L} & 0\\
	0 & 0 & 0 & 0 & 0 & \dfrac{1}{\mu_T}
	\end{pmatrix}
	\begin{pmatrix}
	\sigma_{11} \\ \sigma_{22} \\ \sigma_{33} \\ \sigma_{23} \\ \sigma_{13} \\ \sigma_{12}
	\end{pmatrix}.
	\end{equation}
	In the remainder of this work we make the assumption, with little loss in generality, that $\nu_T = \nu_L := \nu$ and $\mu_T = \mu_L := \mu$. Further, we set 
	\begin{equation}
	p = \frac{E_L}{E_T},
	\label{p}
	\end{equation}
	so that the parameter $p$ measures the degree of transverse isotropy, with inextensible behaviour corresponding to the limit $p \rightarrow \infty$. 
	
	The parameters in \eqref{stress-strain2} may then be expressed in terms of the engineering parameters as \cite{Rasolofoson-Grieshaber-Reddy2018}
	\begin{align}
	\frac{\lambda}{E_T} &=  \dfrac{\nu (\nu + p)}{D},
	\nonumber  \\ \medskip
	\frac{\alpha}{E_T}  &= \dfrac{\nu^2 (p-1)}{D},
	\label{normalitoeng}\\ \medskip
	\frac{\beta}{E_T}   &= \dfrac{p^2(1 - \nu^2) - p(1 + 2\nu^2) + 3\nu^2}{D},
	\nonumber
	\end{align}
	in which the denominator $D$ is given by
	\begin{equation}
	D = (1+\nu)(p(1-\nu) - 2\nu^2)\,.
	\label{denom}
	\end{equation}
	We also have the relation
	\begin{equation}
	\mu_T = \frac{E_T}{2(1 + \nu)}\,.
	\label{muT}
	\end{equation}
	Furthermore, noting that $D \rightarrow 0$ in the incompressible limit $\nu \rightarrow \half$ when $p \rightarrow 1$, which corresponds to the isotropic limit, it is evident from \eqref{normalitoeng} and \eqref{denom} that 
	\begin{align}
	&  \left\{ \begin{array}{l} \mbox{$\lambda$ is bounded as $\nu \rightarrow \half$, if $p > 1$, and as $p \rightarrow \infty$ (inextensibility)} \\ \mbox{$\lambda \rightarrow \infty$ as $\nu \rightarrow \half$, for $p = 1$ (isotropy)}
	\end{array} \right. \\ & \nonumber \\ 
	&  \left\{ \begin{array}{l} \mbox{$\alpha$ is bounded as $\nu \rightarrow \half$, if $p > 1$} \\ \mbox{$\alpha  \rightarrow 0$ as $p \rightarrow 1$ (isotropy)}\end{array} \right. \label{bounds} \\ & \nonumber \\
	&  \left\{ \begin{array}{l} \mbox{$\beta$ is bounded as $\nu \rightarrow \half$, if $p > 1$} \\ \mbox{$\beta  \rightarrow \infty$ as $p \rightarrow \infty$ (inextensibility)}\,.\end{array} \right. 
	\label{limits}
	\end{align}
	
	The elasticity tensor is assumed to be pointwise stable; that is, to satisfy the condition
	\[
	\bepsilon:\BBC\bepsilon > 0 \quad \mbox{for all}\ \bepsilon \,.
	\]
	General conditions on the material constants for pointwise stability are somewhat complex \cite{Rasolofoson-Grieshaber-Reddy2018}, but the simple set of 
	conditions   
	\begin{equation}
	\lambda + \frac{2}{3}\mu > 0,\quad \mu > 0,\quad p \geq 1
	\label{PSconstraints}
	\end{equation}
	meet these requirements. We henceforth assume \eqref{PSconstraints} to hold.
	
	\subsection{Governing equations}
	
	The body is subjected to a body force $\fb$, prescribed loading $\bh$ on $\Gamma_N$, and a prescribed displacement $\bg$ on $\Gamma_D$. 
	
	The equation of equilibrium is 
	\begin{equation}
	-\mbox{div}\,\bsigma = \fb\qquad \mbox{on}\ \Omega.
	\label{equil}
	\end{equation}
	Small displacements are assumed, and the strain displacement relation is
	\begin{equation}
	\bepsilon (\bu) = \half (\nabla\bu + [\nabla\bu]^T)\qquad\mbox{or}\qquad \varepsilon_{ij}(\bu) = \half (u_{i,j} + u_{j,i}).
	\label{strain-disp}
	\end{equation}
	Here $\bu$ denotes the displacement, and $\nabla \bu$ the displacement gradient with components $u_{i,j}$. Here and henceforth we choose a fixed Cartesian coordinate system $x_i$ with orthonormal basis $\be_i$. 
	
	The boundary conditions are
	\begin{subequations}
		\begin{align}
		\bu & = \bg\qquad \mbox{on}\ \Gamma_D, 
		\label{Dirichlet}\\
		\bsigma \cdot \bn & = \bh\qquad \mbox{on}\ \Gamma_N.
		\label{Neumann}
		\end{align}
		\label{bcs}
	\end{subequations}
	Equations \eqref{equil} -- \eqref{bcs}, together with the elastic relation \eqref{stress-strain2}, constitute the boundary-value problem for a transversely isotropic body.
	
	\subsection{Weak formulation}
	We denote by $L^2(\Omega)$ the space of square-integrable functions on $\Omega$, and by $H^1(\Omega)$ the Sobolev space of functions which, together with their generalized first derivatives, are square-integrable, and set $V = [H^1_D(\Omega)]^d = \{ \bv\ |\ v_i \in H^1(\Omega),\ \ \bv = \bzero\ \mbox{on}\  \Gamma_D\}$.
	
	We also define the function $\bu_g \in [ H^1(\Omega)]^d$ such that 
	\[
	\left. \bu_g\right|_{\Gamma_D} = \bg\,.
	\]
	The bilinear form $a(\cdot,\cdot)$ and linear functional $\ell (\cdot)$ are defined by
	\begin{subequations}
		\begin{align}
		a: [ H^1(\Omega)]^d \times [ H^1(\Omega)]^d \rightarrow \BBR,& \hspace{1cm} a(\bu,\bv) = \int_\Omega \bsigma(\bu):\bepsilon (\bv)\ dx,\label{adef}\\
		l: [ H^1(\Omega)]^d \rightarrow \BBR,& \hspace{1cm} \ell (\bv) = \int_\Omega \fb\cdot\bv dx + \int_{\Gamma_N}\bh\cdot\bv\ ds - a(\bu_g,\bv).\label{ldef}
		\end{align}
	\end{subequations}
	The weak form of the problem is then as follows: given $\fb \in [ L^2(\Omega)]^d$ and $\bh \in [ L^2(\Gamma_N)]^d$,
	find $\bU \in [ H^1(\Omega)]^d$ such that $\bU=\bu+\bu_g, \bu \in V$, and
	\begin{equation}
	\label{weak_form}
	a(\bu,\bv) = \ell (\bv) \hspace{1cm} \forall \bv\in V.
	\end{equation}
	We write the bilinear form as
	\begin{equation}
	a(\bu,\bv) = a^{\rm iso}(\bu,\bv) + a^{\rm ti}(\bu,\bv),
	\label{iso+ti}
	\end{equation}
	where
	\begin{subequations}
		\begin{align}
		a^{\rm iso}(\bu,\bv)
		&= \lambda \int_\Omega (\nabla \cdot \bu)(\nabla \cdot \bv) dx
		+ 2 \mu \int_\Omega \bepsilon (\bu):\bepsilon (\bv)\ dx ,\label{a_iso}\\
		a^{\rm ti}(\bu,\bv)
		&=\alpha \int_\Omega \left[(\bM:\bepsilon (\bu))(\nabla \cdot \bv) + (\nabla \cdot \bu)(\bM:\bepsilon (\bv))\right]\ dx
		+ \beta \int_\Omega (\bM:\bepsilon (\bu))(\bM:\bepsilon (\bv))\ dx\,. \nonumber\\
		\label{a_TI}
		\end{align}
		\label{a}
	\end{subequations}
	\vspace*{-5mm}
	\newline \noindent
	The bilinear form is clearly symmetric. The well-posedness of the weak problem requires the bilinear form to be continuous and coercive, and the linear functional continuous. With the assumptions \eqref{PSconstraints}, it is shown in \cite{Rasolofoson-Grieshaber-Reddy2018} that the problem has a unique solution that depends continuously on the data.
	
	\section{The virtual element method}
\begin{figure}[h!]
		\centering
		\includegraphics[scale = 0.33]{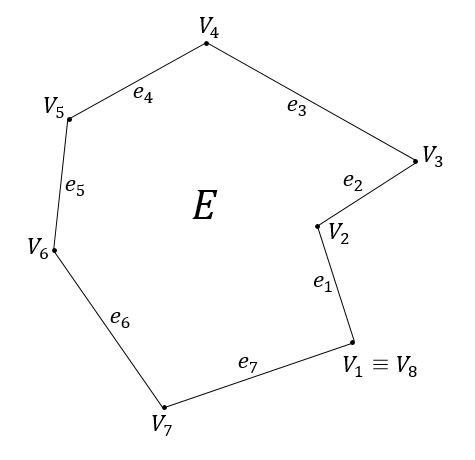}
		\caption{An arbitrary virtual element}
		\label{fig:ArbElement}
	\end{figure}

	The domain $\Omega$ is partitioned into a mesh of elements comprising non-overlapping polygons $E$ with $\overline{\cup E} = \bar{\Omega}$. A typical polygonal element is shown in Figure \ref{fig:ArbElement}.
	A generic edge in a polygonal element is denoted by $e$ and a vertex by $i$, so that $e = 1,\ldots, N$ and $i = 1,\ldots , N$, where $N$ is the number of vertices of an element.
	
	We construct a conforming approximation in a space $V^h \subset V$. The space $V^h$ comprises functions that are continuous on $\Omega$, piecewise linear on the boundary $\partial E$ of each element, and with $\mbox{div}\,\BBC\bepsilon (\bv_h)$ vanishing on $E$ \cite{Beirao-etal2013a,Beirao-etal2014}:
	\begin{equation}
	V^h = \{ \bv_h \in V\ |\ \bv_h \in [C(\Omega)]^2,\ \  \mbox{div}\,\BBC\bepsilon (\bv_h) = \bzero\ \mbox{on}\ E,\ \bv_h|_e \in P_1 (e) \}\,.
	\label{Vh}
	\end{equation}
	Here and henceforth $P_k(X)$ denotes the space of polynomials of degree $\leq k$ on the set $X \subset \BBR^d\ (d =1,2)$. We assign degrees of freedom to the nodes, which are located at the element vertices, and write, for each element,
	\begin{equation}
	\bv_h|_E = \bvarphi \bd
	\label{vh}
	\end{equation}
	in which $\bvarphi$ denotes a matrix of virtual basis functions, and $\bd$ is the $2N \times 1$ vector of degrees of freedom.
	
	All computations will be carried out on the edges $e$ of elements, and it is convenient to write also
	\begin{equation}
	\bv_h|_{\partial E} = \bN \bd\quad \mbox{and}\quad \bepsilon (\bv_h) = \bB\bd\,,
	\label{vhedge}
	\end{equation}
	in which $\bN$ and $\bB$ are respectively matrices of standard Lagrangian linear basis functions and their derivatives. Thus, the basis functions $\bvarphi$ are not known, and are not required to be known; their traces on the boundary are however required, and are simple Lagrangian functions. 
	
	We will require the projection $\Pi: V_h|_E \rightarrow P_0(E)$, defined on $E$ by
	\begin{equation}
	\int_E \Pi \bv_h\ dx = \int_E \bepsilon(\bv_h)\ dx\,.
	\label{Pi}
	\end{equation}
	Thus $\Pi$ is the $L^2$-orthogonal projection onto constants of the strain associated with the displacement $\bv_h$ on an element $E$. 
	From \eqref{vhedge}, and given that $\Pi\bv_h$ is constant we have, in component form,
	\begin{align}
	(\Pi\bv_h)_{ij} & = \frac{1}{2}\frac{1}{|E|} \int_E ((v_h)_{i,j} + (v_h)_{j,i}) \ dx \nonumber \\
	& =  \frac{1}{2}\frac{1}{|E|}  \oint_{\partial E} ((v_h)_in_j + (v_h)_jn_i) \ ds \nonumber \\
	& =  \frac{1}{2}\frac{1}{|E|}\sum_{e \in \partial E} \int_e (N_{iA}d^E_An_j + N_{jA}d^E_An_i)\ ds\,. 
	\label{Pivh}
	\end{align}
	Here $d_A^E$ denotes the degrees of freedom associated with element $E$, summation is implied over all repeated indices, and we have used integration by parts and the representation \eqref{vhedge}$_1$. The integrals in \eqref{Pivh} are readily evaluated as the edge basis functions are known. Thus the projection $\Pi\bv_h$ is available as a function of the degrees of freedom.
	
	To construct the virtual element formulation we start by writing 
	\begin{align}
	a^E (\bu,\bv) & := a(\bv,\bv)|_E \nonumber \\
	& = \int_E \bepsilon(\bu_h) :\BBC\bepsilon(\bv_h)\ dx\,, 
	\label{aE}
	\end{align}
	so that $a^E(\cdot,\cdot)$ denotes the contribution of element $E$ to the bilinear form $a(\cdot,\cdot)$ defined in \eqref{weak_form} and \eqref{a}.
	We have
	\begin{subequations}
		\begin{align}
		a^E(\bu_h,\bv_h) & = \int_E \Pi\bu_h :\BBC\Pi\bv_h\ dx + \int_E (\bepsilon (\bu_h) - \Pi\bu_h) :\BBC(\bepsilon(\bv_h) - \Pi\bv_h)\ dx \label{a} \\
		+ & \int_E \Pi\bu_h :\BBC(\bepsilon(\bv_h) - \Pi\bv_h)\ dx +  \int_E (\bepsilon (\bu_h) - \Pi\bu_h) :\BBC\Pi\bv_h\  dx \label{b} \\
		& = \int_E \Pi\bu_h :\BBC\Pi\bv_h\ dx + \int_E (\bepsilon (\bu_h) - \Pi\bu_h) :\BBC(\bepsilon(\bv_h) - \Pi\bv_h)\ dx \label{c} \\
		& = \underbrace{\int_E \Pi\bu_h :\BBC\Pi\bv_h\ dx}_{\rm consistency\ term} + \underbrace{\int_E \big[ \bepsilon (\bu_h):\BBC\bepsilon(\bv_h) - \Pi\bu_h:\BBC\Pi\bv_h\big]\ dx}_{\rm stabilisation\ term}
		\label{aE2}
		\end{align}
	\end{subequations}
	The last line is obtained by noting the definition of the projection operator, so that the two terms in \eqref{b} are zero. Furthermore, the definition of the projection is invoked again in going from \eqref{c} to \eqref{aE2}. The terms in \eqref{aE2} are referred to respectively as the consistency term and stabilisation term.
	
	\paragraph{The consistency term.}\ \  After substitution of \eqref{Pivh} in the consistency term, evaluation of the integral leads to the expression
	\begin{equation}
	\int_E \Pi\bu_h :\BBC\Pi\bv_h\ dx = (\bar{\bd})^T\bK_{\rm con}^E\bd
	\label{Kcon}
	\end{equation}
	in which $\bK^E_{\rm con}$ is the consistency stiffness matrix for element $E$ and $\bd^E$ and $\bar{\bd}^E$ are respectively the vectors of nodal degrees of freedom of $\bu_h$ and $\bv_h$ on element $E$.

	\paragraph{The stabilisation term.}\ \ Use of the consistency term alone would lead to a rank-deficient stiffness matrix. The second term on the right hand side of \eqref{aE2} serves the purpose of stabilizing the formulation. The basic idea behind the VEM is that integrals are evaluated on the boundaries of elements only; the stabilisation term in its original form would require that integrals be evaluated on the elements. Nevertheless, it is not necessary to evaluate this term exactly, and it suffices to replace it with an approximation. There are several methods that can be employed, see for example \cite{Beirao-etal2015,Gain-etal2014}. However, we choose the stabilisation method presented in \cite{Artioli-etal2017a} as it has proven very robust,  
	\begin{equation}
	a^{E}_{\rm stab} (\bu_h,\bv_h) = \uptau \overline{\bd}^T [\bI - \bD (\bD^T\bD)^{-1}\bD^T ]\bd\,. 
	\label{stabterm}
	\end{equation}
	Here $\overline{\bd}$ and $\bd$ are, again, respectively the vectors of nodal degrees of freedom associated with $\bv_h$ and $\bu_h$, and $\bD$ is the matrix that relates the nodal degrees of freedom 
	$\bd_1$ of a linear vector polynomial to its degrees of freedom $\bs$ relative to a scaled linear monomial basis. That is, for an element with $N$ nodes,
	\begin{equation}
	\bd_1 = {\bD} \bs\,.
	\label{pbases}
	\end{equation}
	Note that $\bD$ has dimensions $2N\times 6$, and has the basis monomials
	\begin{equation}
	\bcM = \left\{ 1, \xi, \eta \right\} = \left\{ 1, \frac{x-x_{c}}{d_{E}},\frac{y-y_{c}}{d_{E}} \right\},
	\end{equation}
	where $d_{E}$ is the diameter of element $E$, with $x_{c}$ and $y_{c}$ the $x-$ and $y-$coordinates of the centroid of $E$ respectively. \\
	This approximation may be motivated by seeking a stabilisation term of the form
	\begin{equation}
	\uptau (\bd^T\bd - \bd_1^T\bd_1)
	\label{newstab}
	\end{equation}
	in which $\bd_1$ are the nodal degrees of freedom of a linear polynomial that is closest to $\bu_h$ in some sense, and $\uptau$ is a suitable scaler to be chosen. 
	In the event that $\bu_h$ is a linear polynomial, then this term vanishes of course.
	
	From \eqref{pbases} we have 
	\begin{align}
	\bd_1^T\bd_1 & = (\bs^T\bD^T)(\bD\bs) &  \nonumber \\
	&  = \bs^T(\bD^T\bD)(\bD^T\bD)^{-1}(\bD^T\bD)\bs & \nonumber \\
	&  = \bd_1^T\bD(\bD^T\bD)^{-1}\bD^T\bd_1\,. 
	\label{d1term}
	\end{align}
	Then we obtain \eqref{stabterm} by replacing $\bd_1$ with the actual vector degrees of freedom.
	
	We need to choose a suitable value for the scalar $\uptau$, such that it is some value representative of the constitutive tensor. We consider the transversely isotropic material properties $\lambda$, $\alpha$, $\beta$, $\mu_{L}$ and $\mu_{T}$. As seen in Section \ref{subsec:ElasticRelation} $\lambda, \beta \rightarrow \infty$ as $\nu \rightarrow0.5$, to keep the VEM locking free we therefore reject these options. We choose $\uptau=\mu_{T}$ as it is bounded and is representative of both isotropic and transversely isotropic materials. As we have set $\mu:=\mu_{L}=\mu_{T}$, we then have
	\begin{equation}
	\boldsymbol{K}_{\rm stab}^{E} = \mu \left[ \boldsymbol{I} - \boldsymbol{\mathcal{D}}\left( \boldsymbol{\mathcal{D}}^{T}\boldsymbol{\mathcal{D}} \right)^{-1}\boldsymbol{\mathcal{D}}^{T} \right]. \label{eqn:StabTerm}
	\end{equation}
	As we have used scaled coordinates, no area scaling of the stabilisation term is necessary.
	The complete stiffness matrix is then given by
	\begin{equation}
	\boldsymbol{K}^{E} = \boldsymbol{K}^{E}_{\rm con} + \boldsymbol{K}^{E}_{\rm stab}.
	\end{equation}
	
	\section{Numerical results}
	\label{sec:NumRes}
	In this section we present numerical results for three model problems to illustrate the performance of the VEM. We consider homogeneous materials, for which the plane of isotropy is fixed across the domain, and also non-homogenous materials, for which the plane of isotropy, as defined by the vector $\ba$, varies with position. 
	Plane strain conditions are assumed. As in Section \ref{subsec:ElasticRelation} we set $\nu_{T}=\nu_{L} = \nu$ and $\mu_{T}=\mu_{L} = \mu$. We consider values of $ p > 1 $ and of Poisson's ratio $\nu = 0.3$ or, to test behaviour in the near-incompressible limit, $\nu = 0.49995$. In all cases the conditions for pointwise stability \eqref{PSconstraints} are met.
	
\begin{figure}[h]
		\centering
		\begin{subfigure}{0.5\textwidth}
			\centering
			\includegraphics[width=0.8\textwidth]{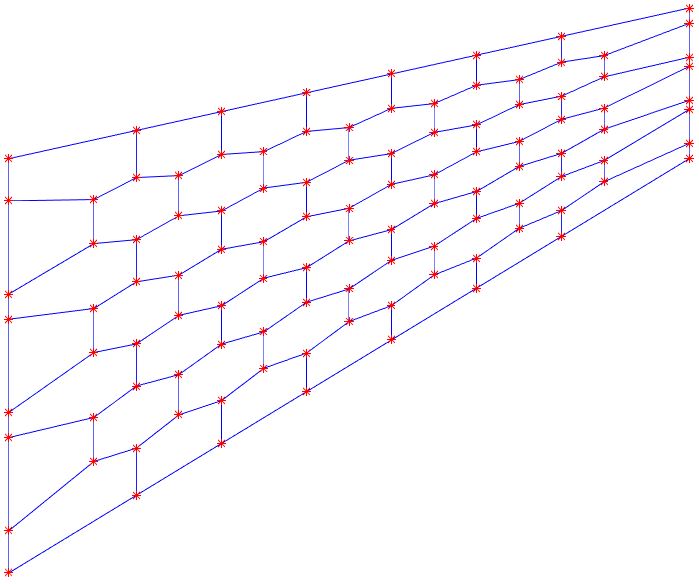}
			\caption{Hex mesh}
			\label{fig:CookHex}
		\end{subfigure}%
		\begin{subfigure}{0.5\textwidth}
			\centering
			\includegraphics[width=0.8\textwidth]{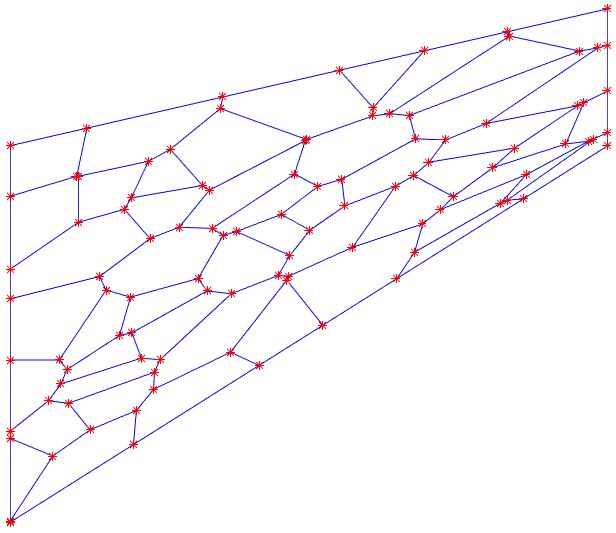}
			\caption{Voronoi mesh}
			\label{fig:CookVoronoi}
		\end{subfigure}
		\caption{Cook's membrane problem, showing the hexagonal and Voronoi meshes for mesh density 7}
		\label{VEM_Mesh_Fig}
	\end{figure}

	We define $\hat{a} := \widehat{(Ox,\boldsymbol{a})}$ to be the angle between the $x$-axis and the fibre direction $\boldsymbol{a}$. The results in the examples that follow are obtained for the following element types:\\
	\vspace*{-5mm}
	\begin{table}[h!]
		\begin{tabular}{lp{15cm}}
			Q$_1$ & The standard bilinear quadrilateral\\
			Q$_2$ & The standard biquadratic quadrilateral\\
			Quad  & The VEM formulation with four-noded elements \\
			Hex & The VEM formulation with six-noded elements \\
			Voronoi & The VEM formulation with Voronoi elements 
		\end{tabular}
	\end{table}
	\vspace*{-5mm}
	\newline \noindent
	Figure \ref{VEM_Mesh_Fig} depicts patches of the meshes comprising six-noded and Voronoi elements for a mesh density $d$ of 7, where $d=\sqrt{n_{elements}}$. Meshes are constructed on a unit domain and then mapped to the problem domain, Figure \ref{VEM_Mesh_Fig} depicts meshes after this mapping.
	
	\medskip
	
	\subsection{Constant fibre direction}
	We present results here for the case in which fibre directions are constant on the domain. The emphasis is on near-incompressibility and near-inextensibility, either separately or combined. In all examples Poisson's ratio is set at $\nu = 0.49995$, and a range of values of $p > 1$ are considered. 
	
	\bigskip
	
\begin{figure}[!htb]
		\centering
		\includegraphics[width=0.4\linewidth]{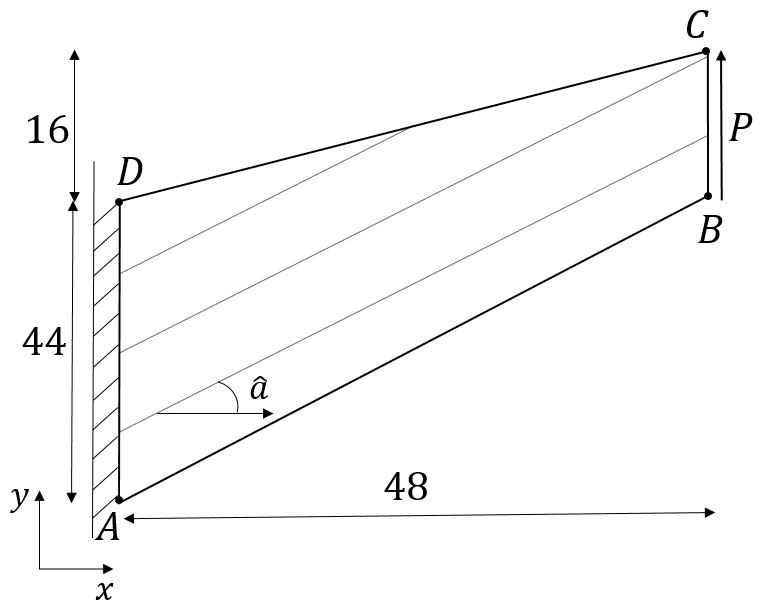}
		\setlength{\belowcaptionskip}{-10pt}
		\caption{Cook's membrane problem, showing fibres inclined at $\hat{a} = \frac{\pi}{4}$}
		\label{Cook_Fig_pi4}
	\end{figure}

	{\bf Cook's membrane problem.} This problem consists of a trapezoidal panel fully fixed along its left edge with a uniformly distributed load along its right edge, as shown in Figure \ref{Cook_Fig_pi4}. The applied load is $ P=100 \rm N $ and  $ E_{T}=250 \rm Pa $. This test problem has no analytical solution. The vertical displacement at  point $ C $ is recorded.
	
	
\begin{figure}[h]
		\centering
		\begin{tikzpicture}
		\begin{axis}[
		xlabel = {Mesh Density},
		ylabel = {Displacement},
		every axis plot/.append style={thick},
		minor y tick num = 1,]
		\addplot[magenta,mark=square*,mark repeat=3] table [x={Density}, y={VEM_quad}] {1_1_1.dat};
		\addlegendentry{Quad};
		\addplot[green,mark=pentagon*,mark repeat=3] table [x={Density}, y={VEM_hex}] {1_1_1.dat};
		\addlegendentry{Hex}
		\addplot[cyan,mark=triangle*,mark repeat=3] table [x={Density}, y={Voronoi}] {1_1_1.dat};
		\addlegendentry{Voronoi}
		\addplot[blue, solid,mark=diamond*,mark repeat=3] table [x={Density}, y={Q1}] {1_1_1.dat};
		\addlegendentry{$Q_{1}$}
		\addplot[black, solid,mark=star, mark repeat=3] table [x={Density}, y={Q2}] {1_1_1.dat};
		\addlegendentry{$Q_{2}$}
		\pgfplotsset{every axis legend/.append style={
				at={(0.97,0.04)},
				anchor=south east}}
		\pgfplotsset{every axis/.append style={
				label style={font=\footnotesize},
				tick label style={font=\footnotesize}  
		}}
		\pgfplotsset{compat=1.3}	
		\end{axis}
		\end{tikzpicture}
		\caption{The Cook problem: convergence test for fibre angle $\hat{a} = \frac{\pi}{4}$ and $p=5$}
		\label{Fig:Cook_Conv_Const}
	\end{figure}
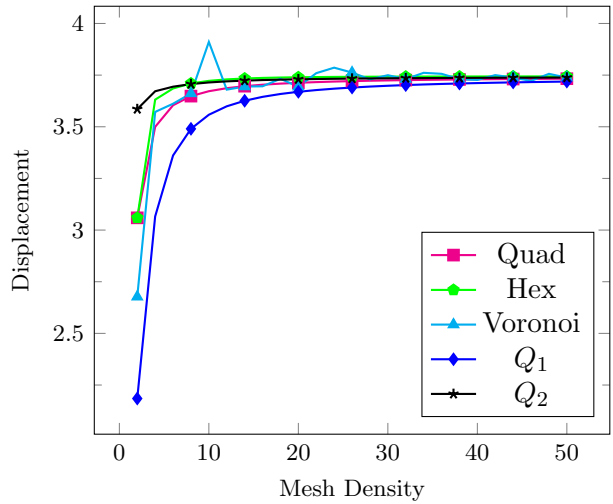

	\noindent Figure \ref{Fig:Cook_Conv_Const} shows a convergence plot of tip displacement vs mesh density for fibre angle $ \hat{a} = \frac{\pi}{4} $, as illustrated in Figure  
	\ref{Cook_Fig_pi4}, and with $p=5$, for the VEM formulation with the three candidate meshes, and for standard finite element formulations. The various VEM formulations are seen to exhibit degrees of accuracy comparable to that of the $Q_2$ approximation. 
	
\begin{figure}[ht!]
		\centering
		\begin{subfigure}{0.5\textwidth}
			\centering
			\begin{tikzpicture}
			\begin{semilogxaxis}[
			xlabel = {$p$},
			ylabel = {Displacement},
			every axis plot/.append style={thick},
			ymin = 0.8, ymax = 8.5,
			minor y tick num = 1,]
			\addplot[magenta,mark=square*] table [x={p}, y={VEM_quad}] {legend.dat};
			\addlegendentry{Quad};
			\addplot[green,mark=pentagon*] table [x={p}, y={VEM_hex}] {legend.dat};
			\addlegendentry{Hex};
			\addplot[cyan,mark=triangle*] table [x={p}, y={Voronoi}] {legend.dat};
			\addlegendentry{Voronoi};
			\addplot[blue,mark=diamond*] table [x={p}, y={Q1}] {legend.dat};
			\addlegendentry{$Q_{1}$};
			\addplot[black,mark=star] table [x={p}, y={Q2}] {legend.dat};
			\addlegendentry{$Q_{2}$};
			\addplot[magenta] table [x={p}, y={VEM_quad}] {1_2_1_1.dat};
			\addplot[only marks,mark=square*,mark options={magenta}] table [x={p}, y={VEM_quad}] {1_2_1_1_m.dat};
			\addplot[green] table [x={p}, y={VEM_hex}] {1_2_1_1.dat};
			\addplot[only marks,mark=pentagon*,mark options={green}] table [x={p}, y={VEM_hex}] {1_2_1_1_m.dat};
			\addplot[cyan] table [x={p}, y={Voronoi}] {1_2_1_1.dat};
			\addplot[only marks,mark=triangle*,mark options={cyan}] table [x={p}, y={Voronoi}] {1_2_1_1_m.dat};
			\addplot[blue, solid] table [x={p}, y={Q1}] {1_2_1_1.dat};
			\addplot[only marks,mark=diamond*,mark options={blue}] table [x={p}, y={Q1}] {1_2_1_1_m.dat};
			\addplot[black, solid] table [x={p}, y={Q2}] {1_2_1_1.dat};
			\addplot[only marks,mark=star,mark options={black}] table [x={p}, y={Q2}] {1_2_1_1_m.dat};
			\draw[->] (axis cs:31.6228,7.0) -- (axis cs:1.0965,7.6);
			\node[draw,align=left] at (axis cs:630.957,7.0) {Quad, Hex \\Voronoi,  $Q_{2}$};
			\draw[->] (axis cs:100,5.5) -- (axis cs:1.047,7.0);
			\node[draw,align=left] at (axis cs:316.228,5.5) {$Q_{1}$};
			\pgfplotsset{every axis/.append style={
					label style={font=\footnotesize},
					tick label style={font=\footnotesize}  
			}}
			\pgfplotsset{every axis legend/.append style={
					at={(0.97,0.2)},
					anchor=south east}}
			\pgfplotsset{compat=1.3}
			\end{semilogxaxis}
			\end{tikzpicture}
		\end{subfigure}%
		\begin{subfigure}{0.5\textwidth}
			\centering
			\begin{tikzpicture}
			\begin{semilogxaxis}[
			xlabel = {$p$},
			ylabel = {Displacement},
			every axis plot/.append style={thick},
			ymin = 1.9, ymax = 8.5,
			minor y tick num = 1,]
			\addplot[magenta,mark=square*] table [x={p}, y={VEM_quad}] {legend.dat};
			\addlegendentry{Quad};
			\addplot[green,mark=pentagon*] table [x={p}, y={VEM_hex}] {legend.dat};
			\addlegendentry{Hex};
			\addplot[cyan,mark=triangle*] table [x={p}, y={Voronoi}] {legend.dat};
			\addlegendentry{Voronoi};
			\addplot[blue,mark=diamond*] table [x={p}, y={Q1}] {legend.dat};
			\addlegendentry{$Q_{1}$};
			\addplot[black,mark=star] table [x={p}, y={Q2}] {legend.dat};
			\addlegendentry{$Q_{2}$};
			\addplot[magenta] table [x={p}, y={VEM_quad}] {1_2_1_2.dat};
			\addplot[only marks,mark=square*,mark options={magenta}] table [x={p}, y={VEM_quad}] {1_2_1_2_m.dat};
			\addplot[green] table [x={p}, y={VEM_hex}] {1_2_1_2.dat};
			\addplot[only marks,mark=pentagon*,mark options={green}] table [x={p}, y={VEM_hex}] {1_2_1_2_m.dat};
			\addplot[cyan] table [x={p}, y={Voronoi}] {1_2_1_2.dat};
			\addplot[only marks,mark=triangle*,mark options={cyan}] table [x={p}, y={Voronoi}] {1_2_1_2_m.dat};
			\addplot[blue, solid] table [x={p}, y={Q1}] {1_2_1_2.dat};
			\addplot[only marks,mark=diamond*,mark options={blue}] table [x={p}, y={Q1}] {1_2_1_2_m.dat};
			\addplot[black, solid] table [x={p}, y={Q2}] {1_2_1_2.dat};
			\addplot[only marks,mark=star,mark options={black}] table [x={p}, y={Q2}] {1_2_1_2_m.dat};
			\draw[->] (axis cs:31.6228,7.0) -- (axis cs:1.096,7.6);
			\node[draw,align=left] at (axis cs:630.957,7.0) {Quad, Hex \\Voronoi, $Q_{2}$};
			\draw[->] (axis cs:100,5.5) -- (axis cs:1.047,7.0);
			\node[draw,align=left] at (axis cs:316.228,5.5) {$Q_{1}$};
			\pgfplotsset{every axis/.append style={
					label style={font=\footnotesize},
					tick label style={font=\footnotesize}  
			}}
			\pgfplotsset{compat=1.3}
			\pgfplotsset{every axis legend/.append style={
					at={(0.97,0.2)},
					anchor=south east}}
			\end{semilogxaxis}
			\end{tikzpicture}
		\end{subfigure}%
		\caption{The Cook problem: tip displacement vs $p$ for (a) fibre direction $\hat{a}=\frac{\pi}{4}$; (b) fibre direction $\hat{a}=\frac{\pi}{9}$}
		\label{Fig:Cook_DISPvsP_Const}
	\end{figure}
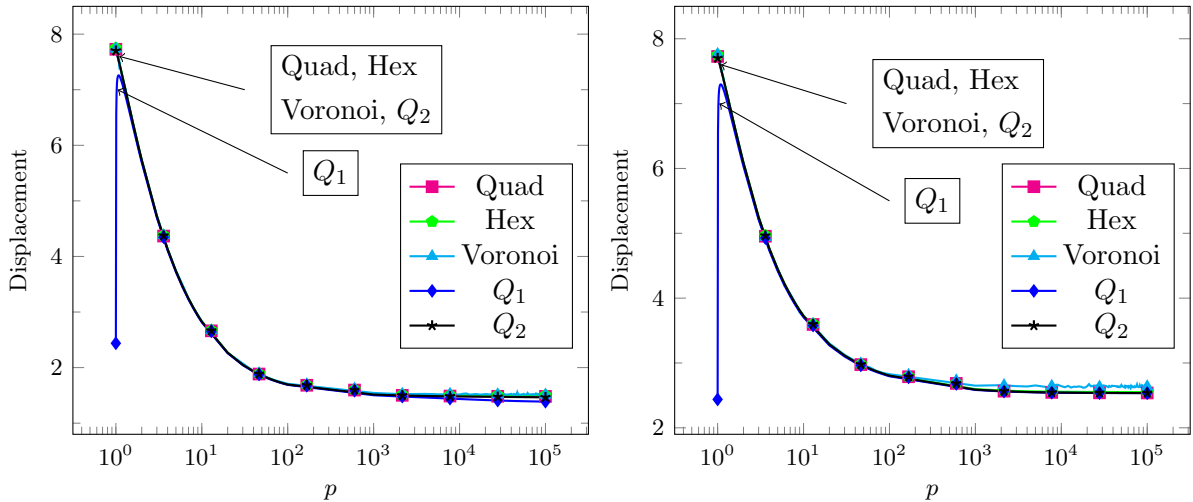

	The results that follow have been generated using meshes with a density of $d=50$.
	\newpage \noindent
	Figure \ref{Fig:Cook_DISPvsP_Const} shows semilog plots of tip displacement vs $p$ for $ 1 \leq p \leq 10^{5} $ and for fibre angles $\hat{a} = \frac{\pi}{4} $ and $\hat{a} = \frac{\pi}{9} $. The well-known locking behaviour of $Q_1$ is clear in the isotropic limit ($ p \rightarrow 1$). On the other hand, the virtual element formulation using quadrilaterals is not equivalent to the conventional formulation using $Q_1$, and is locking-free.
	
\begin{figure}[ht!]
		\centering
		\begin{tikzpicture}
		\begin{axis}[
		xlabel = {$\hat{a}$ (Degrees)},
		ylabel = {Displacement},
		every axis plot/.append style={thick},
		legend columns=1,
		legend pos=outer north east]
		\addplot[magenta,mark=square*,mark repeat=18] table [x={Theta_deg}, y={VEM_quad}] {1_3_1.dat};
		\addlegendentry{Quad}
		\addplot[green,mark=pentagon*,mark repeat=18] table [x={Theta_deg}, y={VEM_hex}] {1_3_1.dat};
		\addlegendentry{Hex}
		\addplot[cyan,mark=triangle*,mark repeat=18] table [x={Theta_deg}, y={Voronoi}] {1_3_1.dat};
		\addlegendentry{Voronoi}
		\addplot[black, solid,mark=star,mark repeat=18] table [x={Theta_deg}, y={Q2}] {1_3_1.dat};
		\addlegendentry{$Q_{2}$}
		\addplot[blue, solid,mark=diamond*,mark repeat=18] table [x={Theta_deg}, y={Q1}] {1_3_1.dat};
		\addlegendentry{$Q_{1}$}
		\draw[->] (axis cs:125,4.0) -- (axis cs:96,3.5);
		\node[draw,align=left] at (axis cs:140,4.0) {$Q_{1}$};
		\pgfplotsset{every axis/.append style={
				label style={font=\footnotesize},
				tick label style={font=\footnotesize}  
		}}
		\pgfplotsset{compat=1.3}	
		\end{axis}
		\end{tikzpicture}
		\caption{The Cook problem: tip displacement vs fibre orientation, for $p= 10^{5}$}
		\label{Fig:Cook_DISPvsTHETA_Const}
	\end{figure}
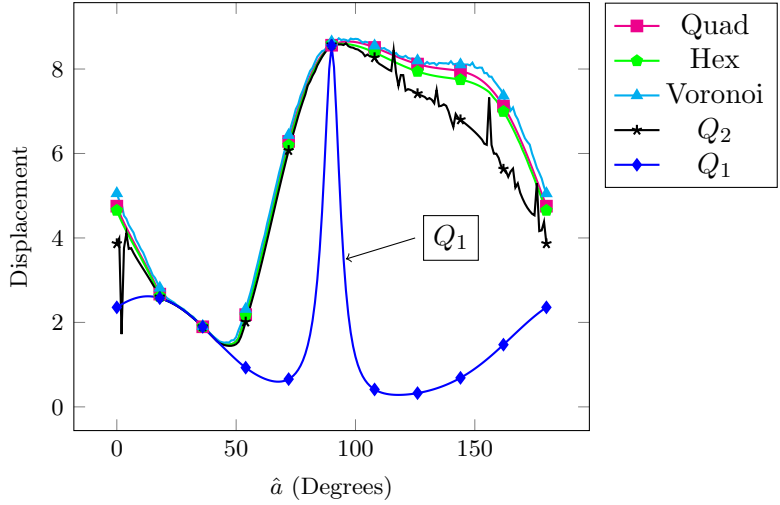

	Figure \ref{Fig:Cook_DISPvsTHETA_Const} shows a plot of tip displacement vs fibre orientation for a nearly inextensible material ($p=10^{5}$). Again, we note the poor performance and locking behaviour of $Q_{1}$ over most of the range, and on the other hand the robust behaviour of the VEM formulation. The $Q_2$ element displays sub-optimal accuracy for fibre angles greater than $\hat{a}=\frac{\pi}{2}$ and close to zero. This is somewhat surprising, in that the behaviour of this element in the near-inextensible limit would be expected to mirror its good performance for near-incompressibility. On the other hand, while the element has been shown to be uniformly convergent for incompressible materials, there does not exist a corresponding analysis for near-inextensibility, to the best of the authors' knowledge. Such an analysis could shed light on the behaviour seen in Figure \ref{Fig:Cook_DISPvsTHETA_Const}.
	\vskip -5mm

	\newpage 
	
\begin{figure}[!h]
		\centering
		\includegraphics[width=0.6\linewidth]{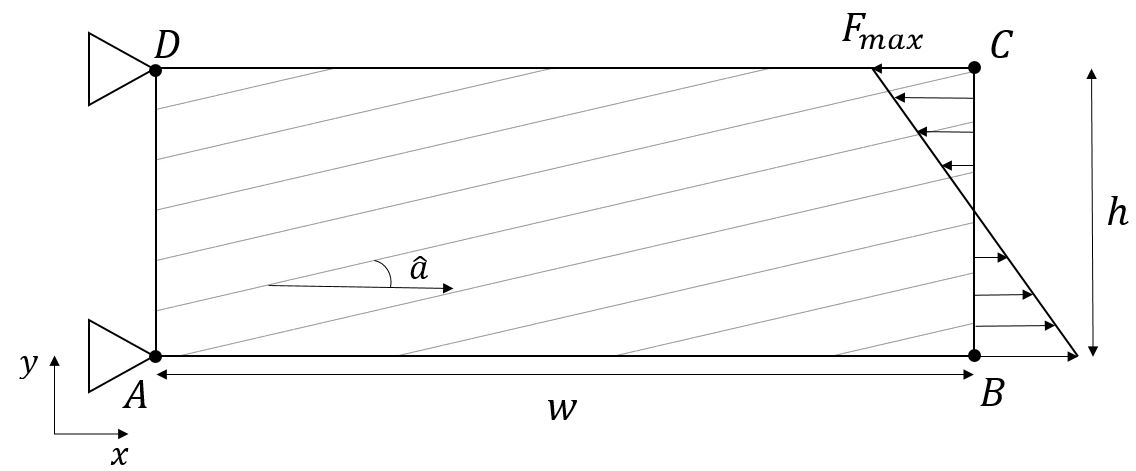}
		\setlength{\belowcaptionskip}{-10pt}
		\caption{The beam problem, showing fibres inclined at $\hat{a} = \frac{\pi}{4}$}
		\label{PureBendA_pi4_Fig}
	\end{figure}

	{\bf The beam problem.} This problem consists of a beam subject to a linearly varying load at its right edge, and pinned at its left extrema, as depicted in Figure \ref{PureBendA_pi4_Fig}. The load has maximum and minimum values of $F_{max}=\pm 30\rm N $. The beam has width $ w=10\rm m $, height $ h=2 \rm m $ and Young's Modulus of $ E_{T} = 1500 \rm Pa $. The vertical displacement at point $ C $ is recorded.
	
	The displacement of point $ C $ is given by \cite{Rasolofoson-Grieshaber-Reddy2018} 
	\begin{align}
	u(x,y) &=\frac{2F_{max}}{h} \left[ \mathbb{S}_{11}xy +\frac{\mathbb{S}_{31}}{2}\left( y^{2}-\frac{h^2}{4} \right) \right] \,, \label{Beam_Solution_X} \\
	v(x,y) &= \frac{F_{max}}{h}\left[ \mathbb{S}_{21} \left( y^{2} - \frac{h^{2}}{4}\right) - \mathbb{S}_{11}x^{2}
	\right]\,; \label{Beam_Solution_Y}
	\end{align}
	the coefficients $\mathbb{S}_{ij}$ are lengthy functions of the material constants, and are given in the Appendix to \cite{Rasolofoson-Grieshaber-Reddy2018}. 
	
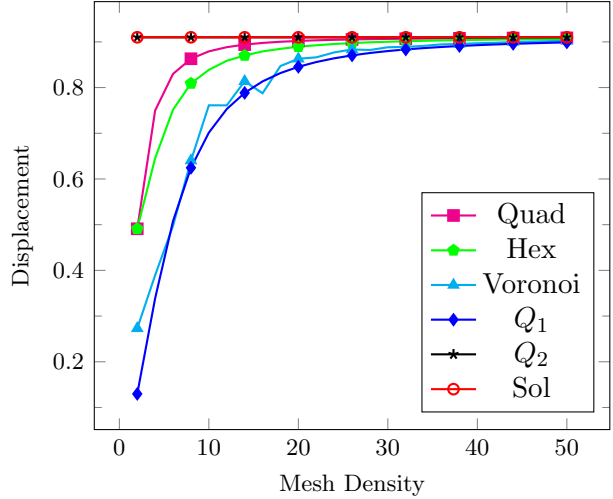
\begin{figure}[h]
		\centering
		\centering
		\begin{tikzpicture}
		\begin{axis}[
		xlabel = {Mesh Density},
		ylabel = {Displacement},
		minor y tick num = 1,
		every axis plot/.append style={thick}]
		\addplot[magenta,mark=square*,mark repeat=3] table [x={Density}, y={VEM_quad}] {1_1_2.dat};
		\addlegendentry{Quad}
		\addplot[green,mark=pentagon*,mark repeat=3] table [x={Density}, y={VEM_hex}] {1_1_2.dat};
		\addlegendentry{Hex}
		\addplot[cyan,mark=triangle*,mark repeat=3] table [x={Density}, y={Voronoi}] {1_1_2.dat};
		\addlegendentry{Voronoi}
		\addplot[blue, solid,mark=diamond*,mark repeat=3] table [x={Density}, y={Q1}] {1_1_2.dat};
		\addlegendentry{$Q_{1}$}
		\addplot[black, solid,mark=star,mark repeat=3] table [x={Density}, y={Q2}] {1_1_2.dat};
		\addlegendentry{$Q_{2}$}
		\addplot[red,mark=o,mark repeat=3] table [x={Density}, y={Analytical}] {1_1_2.dat};
		\addlegendentry{Sol}
		\pgfplotsset{every axis legend/.append style={
				at={(0.97,0.04)},
				anchor=south east}}
		\pgfplotsset{every axis/.append style={
				label style={font=\footnotesize},
				tick label style={font=\footnotesize}  
		}}
		\pgfplotsset{compat=1.3}	
		\end{axis}
		\end{tikzpicture}
		\caption{The beam problem: convergence test for fibre angle $\hat{a}=\frac{\pi}{4}$ and $p=5$}
		\label{Fig:Beam_Conv_Const}
	\end{figure}

	Figure \ref{Fig:Beam_Conv_Const} shows a convergence plot of tip displacement vs mesh density for a fibre orientation of $\hat{a} = \frac{\pi}{4} $, and with $p=5$. It is seen that for the various VEM meshes the convergence behaviour is similar to that of the $Q_2$ mesh for sufficiently fine meshes.
	\newpage
	
	The results that follow have been generated using meshes with a density of $d=50$. 
	
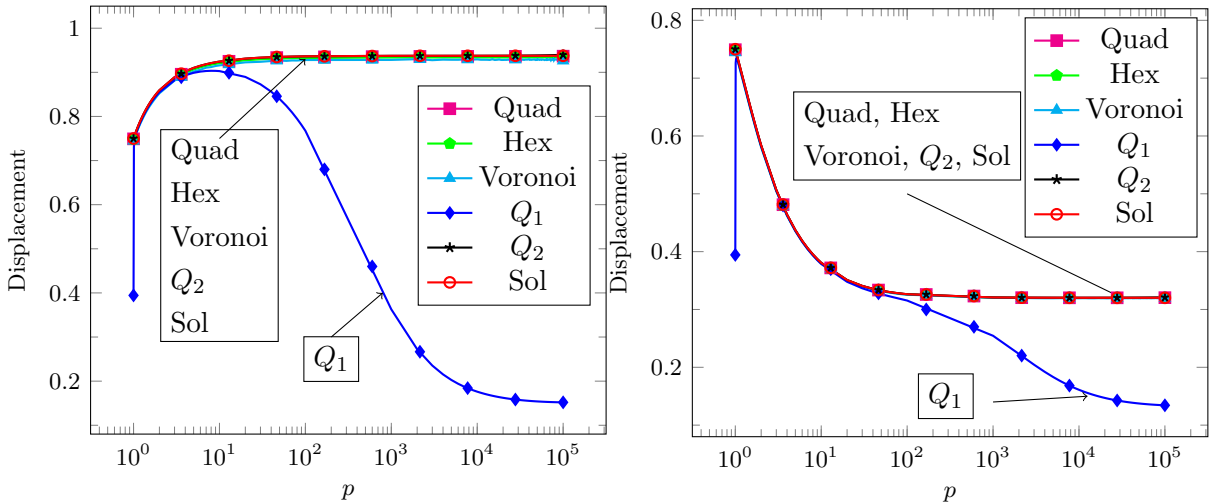
\begin{figure}[ht!]
		\centering
		\begin{subfigure}{0.5\textwidth}
			\centering
			\begin{tikzpicture}
			\begin{semilogxaxis}[
			xlabel = {$p$},
			ylabel = {Displacement},
			every axis plot/.append style={thick},
			ymin = 0.08, ymax = 1.05,
			minor y tick num = 1,]
			\addplot[magenta,mark=square*] table [x={p}, y={VEM_quad}] {legend.dat};
			\addlegendentry{Quad};
			\addplot[green,mark=pentagon*] table [x={p}, y={VEM_hex}] {legend.dat};
			\addlegendentry{Hex};
			\addplot[cyan,mark=triangle*] table [x={p}, y={Voronoi}] {legend.dat};
			\addlegendentry{Voronoi};
			\addplot[blue,mark=diamond*] table [x={p}, y={Q1}] {legend.dat};
			\addlegendentry{$Q_{1}$};
			\addplot[black,mark=star] table [x={p}, y={Q2}] {legend.dat};
			\addlegendentry{$Q_{2}$};
			\addplot[red,mark=o] table [x={p}, y={Sol}] {legend.dat};
			\addlegendentry{Sol};
			\addplot[magenta] table [x={p}, y={VEM_quad}] {1_2_2_1.dat};
			\addplot[only marks,mark=square*,mark options={magenta}] table [x={p}, y={VEM_quad}] {1_2_2_1_m.dat};
			\addplot[green] table [x={p}, y={VEM_hex}] {1_2_2_1.dat};
			\addplot[only marks,mark=pentagon*,mark options={green}] table [x={p}, y={VEM_hex}] {1_2_2_1_m.dat};
			\addplot[cyan] table [x={p}, y={Voronoi}] {1_2_2_1.dat};
			\addplot[only marks,mark=triangle*,mark options={cyan}] table [x={p}, y={Voronoi}] {1_2_2_1_m.dat};
			\addplot[blue, solid] table [x={p}, y={Q1}] {1_2_2_1.dat};
			\addplot[only marks,mark=diamond*,mark options={blue}] table [x={p}, y={Q1}] {1_2_2_1_m.dat};
			\addplot[black, solid] table [x={p}, y={Q2}] {1_2_2_1.dat};
			\addplot[only marks,mark=star,mark options={black}] table [x={p}, y={Q2}] {1_2_2_1_m.dat};
			\addplot[red, solid] table [x={p}, y={Analytical}] {1_2_2_1.dat};
			\addplot[only marks,mark=o,mark options={red}] table [x={p}, y={Sol}] {1_2_2_1_m.dat};
			\draw[->] (axis cs:10,0.75) -- (axis cs:100,0.93);
			\node[draw,align=left] at (axis cs:10,0.53) {Quad\\ Hex \\ Voronoi \\ $Q_{2}$ \\ Sol};
			\draw[->] (axis cs:200,0.3) -- (axis cs:794.328,0.4);
			\node[draw,align=left] at (axis cs:200,0.25) {$Q_{1}$};
			\pgfplotsset{every axis/.append style={
					label style={font=\footnotesize},
					tick label style={font=\footnotesize}  
			}}
			\pgfplotsset{compat=1.3}
			\pgfplotsset{every axis legend/.append style={
					at={(0.97,0.3)},
					anchor=south east}}
			\end{semilogxaxis}
			\end{tikzpicture}
			\label{Fig:Beam_DISPvsP_pi4_Const}
		\end{subfigure}%
		\begin{subfigure}{0.5\textwidth}
			\centering
			\begin{tikzpicture}
			\begin{semilogxaxis}[
			xlabel = {$p$},
			ylabel = {Displacement },
			every axis plot/.append style={thick},
			ymin = 0.08, ymax = 0.82,
			minor y tick num = 1,]
			\addplot[magenta,mark=square*] table [x={p}, y={VEM_quad}] {legend.dat};
			\addlegendentry{Quad};
			\addplot[green,mark=pentagon*] table [x={p}, y={VEM_hex}] {legend.dat};
			\addlegendentry{Hex};
			\addplot[cyan,mark=triangle*] table [x={p}, y={Voronoi}] {legend.dat};
			\addlegendentry{Voronoi};
			\addplot[blue,mark=diamond*] table [x={p}, y={Q1}] {legend.dat};
			\addlegendentry{$Q_{1}$};
			\addplot[black,mark=star] table [x={p}, y={Q2}] {legend.dat};
			\addlegendentry{$Q_{2}$};
			\addplot[red,mark=o] table [x={p}, y={Sol}] {legend.dat};
			\addlegendentry{Sol};
			\addplot[magenta] table [x={p}, y={VEM_quad}] {1_2_2_2.dat};
			\addplot[only marks,mark=square*,mark options={magenta}] table [x={p}, y={VEM_quad}] {1_2_2_2_m.dat};
			\addplot[green] table [x={p}, y={VEM_hex}] {1_2_2_2.dat};
			\addplot[only marks,mark=pentagon*,mark options={green}] table [x={p}, y={VEM_hex}] {1_2_2_2_m.dat};
			\addplot[cyan] table [x={p}, y={Voronoi}] {1_2_2_2.dat};
			\addplot[only marks,mark=triangle*,mark options={cyan}] table [x={p}, y={Voronoi}] {1_2_2_2_m.dat};
			\addplot[blue, solid] table [x={p}, y={Q1}] {1_2_2_2.dat};
			\addplot[only marks,mark=diamond*,mark options={blue}] table [x={p}, y={Q1}] {1_2_2_2_m.dat};
			\addplot[black, solid] table [x={p}, y={Q2}] {1_2_2_2.dat};
			\addplot[only marks,mark=star,mark options={black}] table [x={p}, y={Q2}] {1_2_2_2_m.dat};
			\addplot[red, solid] table [x={p}, y={Analytical}] {1_2_2_2.dat};
			\addplot[only marks,mark=o,mark options={red}] table [x={p}, y={Sol}] {1_2_2_2_m.dat};
			\draw[->] (axis cs:100,0.5) -- (axis cs:31622.78,0.32);
			\node[draw,align=left] at (axis cs:100,0.6) {Quad, Hex \\ Voronoi, $Q_{2}$, Sol};
			\draw[->] (axis cs:1000,0.14) -- (axis cs:12589.25,0.15);
			\node[draw,align=left] at (axis cs:281.828,0.15) {$Q_{1}$};
			\pgfplotsset{every axis/.append style={
					label style={font=\footnotesize},
					tick label style={font=\footnotesize}  
			}}
			\pgfplotsset{compat=1.3}
			\end{semilogxaxis}
			\end{tikzpicture}
			\label{Fig:Beam_DISPvsP_pi9_Const}
		\end{subfigure}%
		\caption{Beam problem: tip displacement vs $p$ for fibre angle (a) $\hat{a}=\frac{\pi}{4}$; (b) $\hat{a}=\frac{\pi}{9}$}
		\label{PureBend_Const}
	\end{figure}

	Figure \ref{PureBend_Const} shows semilog plots of tip displacement vs $p$ for $ 1 \leq p \leq 10^{5} $ for fibre angles of $ \hat{a}=\frac{\pi}{4} $ and $\hat{a} = \frac{\pi}{9} $. Again, as with the Cook problem, the VEM solutions are locking-free and display high accuracy. As pointed out in \cite{Rasolofoson-Grieshaber-Reddy2018}, for mild anisotropy, that is, low values of $p$, the tendency to lock for the $Q_1$ mesh is mitigated as a result of the Lam\'e parameter being bounded for $p>1$ in conditions of near-incompressibility, with $\nu$ very close to 0.5. This behaviour is evident in Figure \ref{PureBend_Const}, where the $Q_1$ mesh is seen to be locking-free for $p>1$  and for $p$ up to $p \approx 10$ for fibre angle $\hat{a}=\frac{\pi}{4}$, and $p\approx 100$ for fibre angle $\hat{a}=\frac{\pi}{9}$.
	
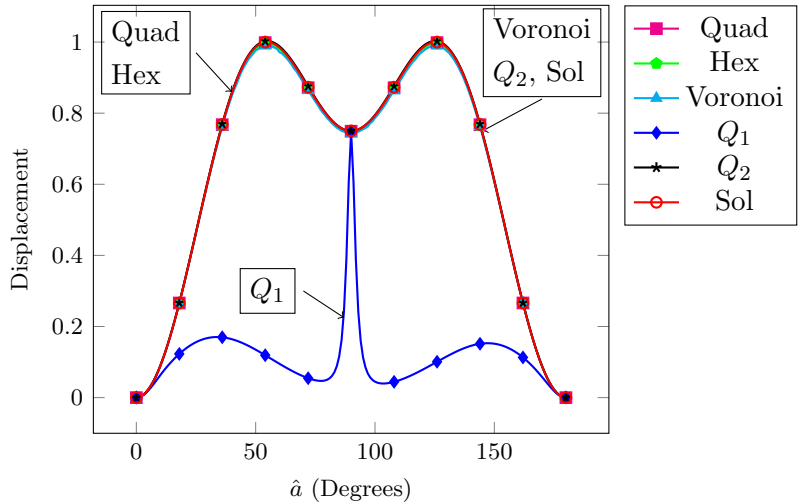
\begin{figure}[h!]
		\centering
		\begin{tikzpicture}
		\begin{axis}[
		xlabel = {$\hat{a}$ (Degrees)},
		ylabel = {Displacement},
		every axis plot/.append style={thick},
		legend columns=1,
		legend pos=outer north east]
		\addplot[magenta,mark=square*,mark repeat=18] table [x={Theta_deg}, y={VEM_quad}] {1_3_2.dat};
		\addlegendentry{Quad};
		\addplot[green,mark=pentagon*,mark repeat=18] table [x={Theta_deg}, y={VEM_hex}] {1_3_2.dat};
		\addlegendentry{Hex};
		\addplot[cyan,mark=triangle*,mark repeat=18] table [x={Theta_deg}, y={Voronoi}] {1_3_2.dat};
		\addlegendentry{Voronoi};
		\addplot[blue, solid,mark=diamond*,mark repeat=18] table [x={Theta_deg}, y={Q1}] {1_3_2.dat};
		\addlegendentry{$Q_{1}$};
		\addplot[black, solid,mark=star,mark repeat=18] table [x={Theta_deg}, y={Q2}] {1_3_2.dat};
		\addlegendentry{$Q_{2}$};
		\addplot[red, solid,mark=o,mark repeat=18] table [x={Theta_deg}, y={Analytical}] {1_3_2.dat};
		\addlegendentry{Sol};
		\draw[->] (axis cs:25,0.97) -- (axis cs:40,0.86);
		\node[draw,align=left] at (axis cs:4.0,0.97) {Quad \\ Hex};
		\draw[->] (axis cs:170,0.84) -- (axis cs:145,0.75);
		\node[draw,align=left] at (axis cs:170.0,0.97) {Voronoi \\ $Q_{2}$, Sol};
		\draw[->] (axis cs:70,0.3) -- (axis cs:87,0.225);
		\node[draw,align=left] at (axis cs:55,0.3) {$Q_{1}$};
		\pgfplotsset{every axis/.append style={
				label style={font=\footnotesize},
				tick label style={font=\footnotesize}  
		}}
		\pgfplotsset{compat=1.3}
		\end{axis}
		\end{tikzpicture}
		\caption{The beam problem: tip displacement vs fibre orientation, for $p= 10^{5}$}
		\label{Fig:Beam_DISPvsTHETA_Const}
	\end{figure}

	Figure \ref{Fig:Beam_DISPvsTHETA_Const} shows a plot of tip displacement vs fibre orientation for the case of near-inextensibility ($p=10^{5}$). Again, we note poor performance of $Q_{1}$ and robust and accurate behaviour of the VEM formulations. In contrast to the result for the Cook problem, here the $Q_2$ element demonstrates equally accurate behaviour.
	
	
	\subsection{Non-homogeneous materials: variable fibre orientation}
\begin{figure}[!h]
		\centering
		\begin{tikzpicture}
		\begin{axis}[
		xlabel = {Mesh Density},
		ylabel = {$w_{c}$},
		minor y tick num = 1,
		every axis plot/.append style={thick}]
		\addplot[blue] table [x={X}, y={Y}] {WeightingFunction.dat};
		\draw [thick,dashed] (axis cs: 10,-0.03) -- (axis cs: 10,0.53);
		\node[text width=3cm] at (axis cs: 14.5,0.01) 
		{$d_{\rm cr}$};
		\pgfplotsset{every axis/.append style={
				label style={font=\footnotesize},
				tick label style={font=\footnotesize}  
		}}
		\pgfplotsset{compat=1.3}
		\end{axis}
		\end{tikzpicture}
		\caption{The weight $w_c$ as a function of mesh density $d$ }
		\label{w_vs_d}
	\end{figure}
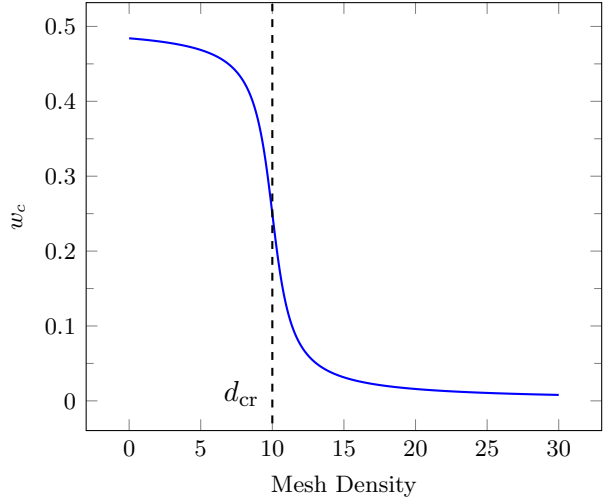

	For a given distribution of fibre directions $\ba (\boldsymbol{\bx})$ it follows that some approximation has to be made for $\ba(\bx)$ within each element so as to preserve the general approach to carrying out VEM computations. A simple option would be to  approximate $\ba$ by its centroidal value. However, such an approximation can be somewhat inaccurate for situations in which the fibre orientation varies significantly across a length scale comparable to mesh size.
	A more reliable approach is to use the average fibre direction at the element nodes; this approach is observed to yield more stable and faster convergence.
	When dealing with rapidly varying fibre directions a more stable approach was achieved by using a weighted average of the fibre direction at the centroid and the average direction at the vertices. This approach applies an equal weighting to the centroidal direction and the average of its nodal values for very coarse meshes; for finer meshes, that is, as the mesh density increases, the weighting of the centroidal value decreases rapidly. Thus for very fine meshes and rapidly varying fibre directions, it is the nodal average that dominates. With the mesh density denoted by $d$, the centroidal weight $w_c$ is defined by  
	\begin{align}
	w_c &= \frac{\frac{\pi}{2} + \arctan(d_{\rm cr} - d)}{2 \pi} \,.
	\label{w_c}
	\end{align}
	This function is shown in Figure \ref{w_vs_d}. Here $d_{\rm cr}$ is a user-defined critical mesh density beyond which the value of the weight drops rapidly. 
	
	The average fibre direction $\ba$ is then given by 
	\begin{equation}
	\ba_{\rm ave}|_E = w_c \ba (\bx_{c}) + (1 - w_c)\frac{1}{N}\sum_{i=1}^{N} \ba (\bx_i)\,, 
	\label{a_ave}
	\end{equation}
	where, as before, $N$ denotes the number of nodes of element $E$, and $\bx_c$ and $\bx_i$ are respectively the coordinates of the element centroid and node $i$.
	We then approximate the elasticity tensor on an element $E$ by 
	\begin{equation}
	\mathbb{C}|_{E} \simeq \mathbb{C}(\boldsymbol{a}_{\rm ave}) \,.
	\label{C_ave}
	\end{equation}
	The critical density used is problem specific as it depends on the degree of variation in fibre orientation. However, for simplicity a critical density of $ {d_{\rm crit}}=10 $ was used as it worked well across a range of problems.
	
	Except where otherwise stated, Poisson's ratio is set at $\nu = 0.3$ in the examples that follow.
	
\begin{figure}[!ht]
		\centering
		\begin{minipage}{.5\textwidth}
			\centering
			\includegraphics[width=0.75\textwidth]{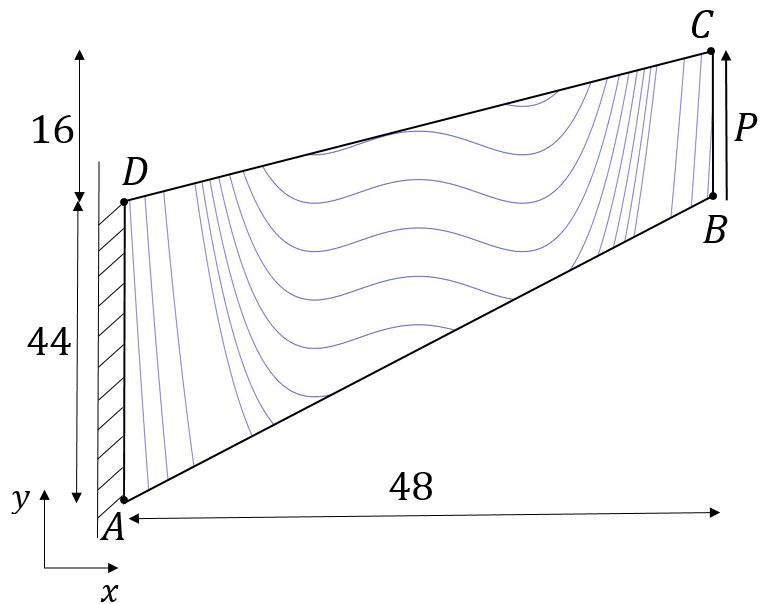}
			\setlength{\belowcaptionskip}{-10pt}
		\end{minipage}%
		\begin{minipage}{0.5\textwidth}
			\centering
			\includegraphics[width=0.75\linewidth]{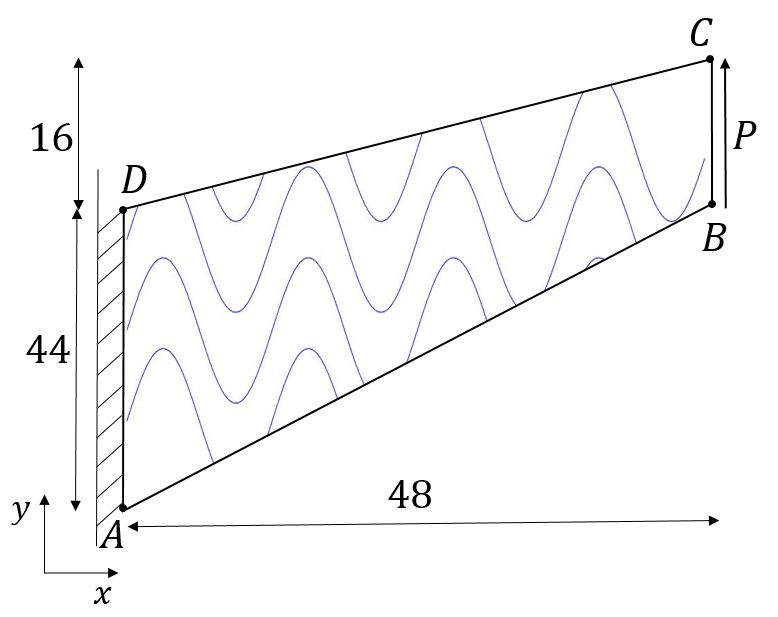}
			\setlength{\belowcaptionskip}{-10pt}
		\end{minipage}
		\caption{Cook's membrane problem with curves showing variable fibre orientation for (a) quartic, and (b) sinusoidal distributions}
		\label{CookFigsNonConst}
	\end{figure}

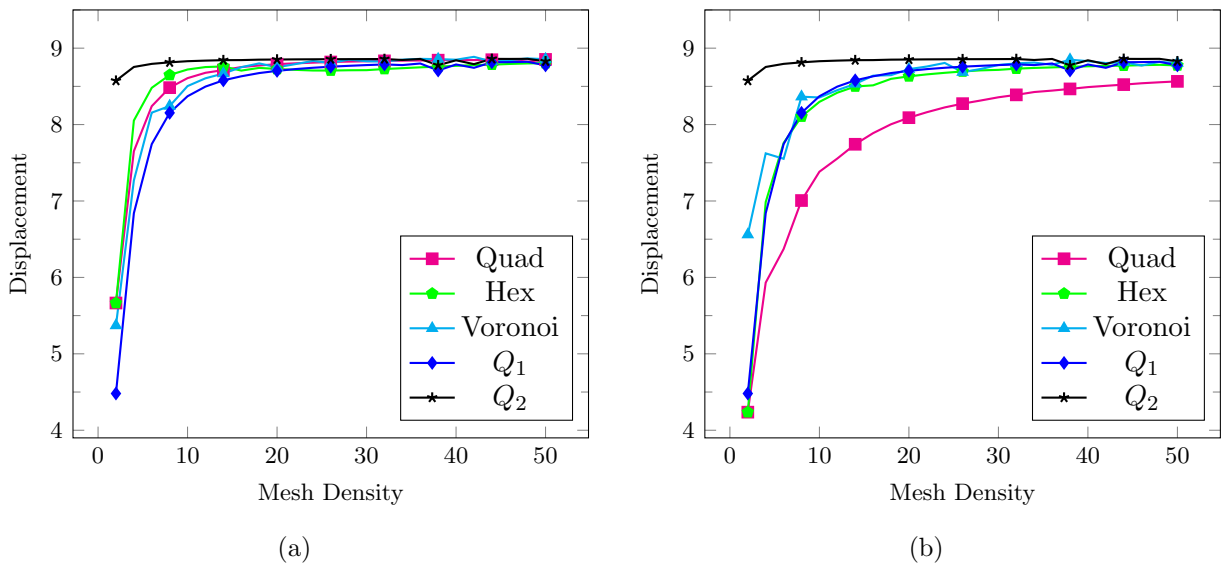
\begin{figure}[h!]
		\centering
		\begin{subfigure}{0.55\textwidth}
			\centering
			\begin{tikzpicture}
			\begin{axis}[
			xlabel = {Mesh Density},
			ylabel = {Displacement},
			ymin = 3.9, ymax = 9.5,
			minor y tick num = 1,
			every axis plot/.append style={thick}]
			\addplot[magenta,mark=square*,mark repeat=3] table [x={Density}, y={VEM_quad}] {Cook_TI_NC_poly4_Cent.dat};
			\addlegendentry{Quad}
			\addplot[green,mark=pentagon*,mark repeat=3] table [x={Density}, y={VEM_hex}] {Cook_TI_NC_poly4_Cent.dat};
			\addlegendentry{Hex}
			\addplot[cyan,mark=triangle*,mark repeat=3] table [x={Density}, y={Voronoi}] {Cook_TI_NC_poly4_Cent.dat};
			\addlegendentry{Voronoi}
			\addplot[blue, solid,mark=diamond*,mark repeat=3] table [x={Density}, y={Q1}] {Cook_TI_NC_poly4_Cent.dat};
			\addlegendentry{$Q_{1}$}
			\addplot[black, solid,mark=star,mark repeat=3] table [x={Density}, y={Q2}] {Cook_TI_NC_poly4_Cent.dat};
			\addlegendentry{$Q_{2}$}
			\pgfplotsset{every axis legend/.append style={
					at={(0.97,0.04)},
					anchor=south east}}
			\pgfplotsset{every axis/.append style={
					label style={font=\footnotesize},
					tick label style={font=\footnotesize}  
			}}
			\pgfplotsset{compat=1.3}	
			\end{axis}
			\end{tikzpicture}
			\caption{}
		\end{subfigure}%
		\begin{subfigure}{0.5\textwidth}
			\centering
			\begin{tikzpicture}
			\begin{axis}[
			xlabel = {Mesh Density},
			ylabel = {Displacement},
			ymin = 3.9, ymax = 9.5,
			minor y tick num = 1,
			every axis plot/.append style={thick}]
			\addplot[magenta,mark=square*,mark repeat=3] table [x={Density}, y={VEM_quad}] {Cook_TI_NC_poly4_Vert.dat};
			\addlegendentry{Quad}
			\addplot[green,mark=pentagon*,mark repeat=3] table [x={Density}, y={VEM_hex}] {Cook_TI_NC_poly4_Vert.dat};
			\addlegendentry{Hex}
			\addplot[cyan,mark=triangle*,mark repeat=3] table [x={Density}, y={Voronoi}] {Cook_TI_NC_poly4_Vert.dat};
			\addlegendentry{Voronoi}
			\addplot[blue, solid,mark=diamond*,mark repeat=3] table [x={Density}, y={Q1}] {Cook_TI_NC_poly4_Vert.dat};
			\addlegendentry{$Q_{1}$}
			\addplot[black, solid,mark=star,mark repeat=3] table [x={Density}, y={Q2}] {Cook_TI_NC_poly4_Vert.dat};
			\addlegendentry{$Q_{2}$}
			\pgfplotsset{every axis legend/.append style={
					at={(0.97,0.04)},
					anchor=south east}}
			\pgfplotsset{every axis/.append style={
					label style={font=\footnotesize},
					tick label style={font=\footnotesize}  
			}}
			\pgfplotsset{compat=1.3}	
			\end{axis}
			\end{tikzpicture}
			\caption{}
		\end{subfigure}
		\setlength{\belowcaptionskip}{-20pt}
		\caption{Tip displacement vs mesh density for Cook's membrane problem for $p=5$, with fibre directions defined by quartic curves, and using (a) the centroidal value of $\ba$; (b) the average of nodal values of $\ba$}
		\label{Fig:CookPoly4}
	\end{figure}

	We present results for two families of fibre distributions, corresponding to curves $y = c + f(x)$ where $f(x)$ is chosen to be, respectively, $(x-24)^{2}(x-12)(x-36)$ and $2\sin x$. The polynomial distribution corresponds to mild variation with position, while the sinusoidal distribution is a more severe test of performance under conditions of rapidly varying direction. Figure \ref{CookFigsNonConst} shows schematically the curves corresponding to these two cases for the Cook problem, one of the examples considered in what follows.
	\vspace*{-3mm}
	\newpage
	{\bf Cook's membrane problem.} Figures \ref{Fig:CookPoly4} (a) and (b) show the tip displacement as a function of mesh density for fibres corresponding to the quartic distribution, with $p=5$, and in which the value of $\ba$ is based respectively on its value at the element centroids and the average of its values at the nodes. We note smooth and stable convergence of the VEM, with the quadrilateral VEM mesh performing somewhat more poorly for the case in which the average nodal value of $\ba$ is used.
	
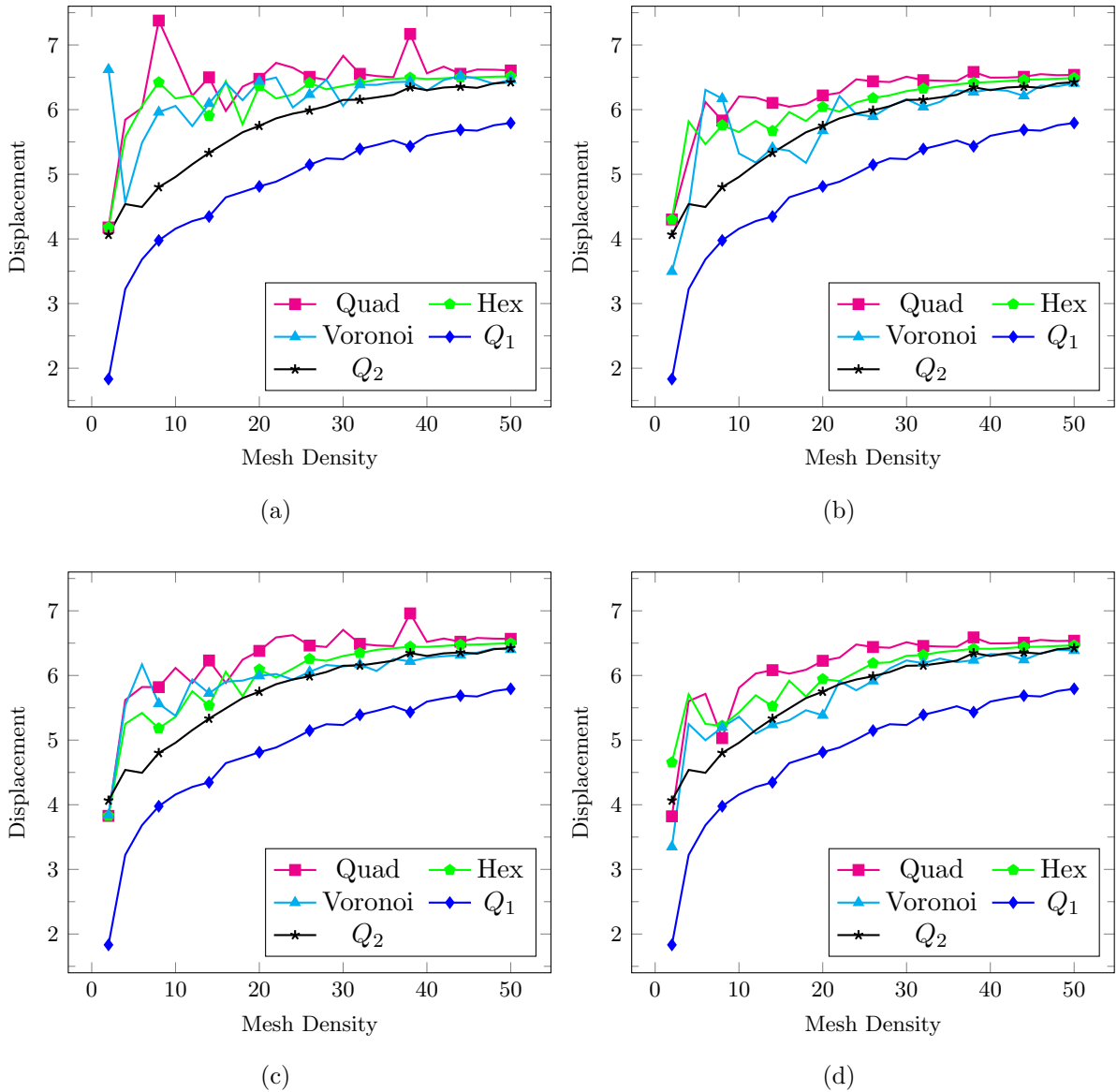
\begin{figure}[ht!]
		\centering
		\begin{subfigure}{0.5\textwidth}
			\centering
			\begin{tikzpicture}
			\begin{axis}[
			xlabel = {Mesh Density},
			ylabel = {Displacement},
			ymin = 1.4, ymax = 7.6,
			minor y tick num = 1,
			ytick = {1, 2, 3, 4, 5, 6, 7},
			every axis plot/.append style={thick},
			legend columns=2, 
			legend style={
				/tikz/column 2/.style={
					column sep=3.5pt}}]
			\pgfplotsset{every axis legend/.append style={	
					at={(0.97,0.04)},
					anchor=south east}}]
			\addplot[magenta,mark=square*,mark repeat=3] table [x={Density}, y={VEM_quad}] {Cook_TI_NC_2sinx_Cent.dat};
			\addlegendentry{Quad}
			\addplot[green,mark=pentagon*,mark repeat=3] table [x={Density}, y={VEM_hex}] {Cook_TI_NC_2sinx_Cent.dat};
			\addlegendentry{Hex}
			\addplot[cyan,mark=triangle*,mark repeat=3] table [x={Density}, y={Voronoi}] {Cook_TI_NC_2sinx_Cent.dat};
			\addlegendentry{Voronoi}
			\addplot[blue, solid,mark=diamond*,mark repeat=3] table [x={Density}, y={Q1}] {Cook_TI_NC_2sinx_Cent.dat};
			\addlegendentry{$Q_{1}$}
			\addplot[black, solid,mark=star,mark repeat=3] table [x={Density}, y={Q2}] {Cook_TI_NC_2sinx_Cent.dat};
			\addlegendentry{$Q_{2}$}
			\pgfplotsset{every axis/.append style={
					label style={font=\footnotesize},
					tick label style={font=\footnotesize}  
			}}
			\pgfplotsset{compat=1.3}
			\end{axis}
			\end{tikzpicture}
			\caption{}
		\end{subfigure}%
		\begin{subfigure}{0.5\textwidth}
			\centering
			\begin{tikzpicture}
			\begin{axis}[
			xlabel = {Mesh Density},
			ylabel = {Displacement},
			ymin = 1.4, ymax = 7.6,
			minor y tick num = 1,
			ytick = {1, 2, 3, 4, 5, 6, 7},
			every axis plot/.append style={thick},
			legend columns=2, 
			legend style={
				/tikz/column 2/.style={
					column sep=3.5pt}}]
			\pgfplotsset{every axis legend/.append style={	
					at={(0.97,0.04)},
					anchor=south east}}
			\addplot[magenta,mark=square*,mark repeat=3] table [x={Density}, y={VEM_quad}] {Cook_TI_NC_2sinx_Vert.dat};
			\addlegendentry{Quad}
			\addplot[green,mark=pentagon*,mark repeat=3] table [x={Density}, y={VEM_hex}] {Cook_TI_NC_2sinx_Vert.dat};
			\addlegendentry{Hex}
			\addplot[cyan,mark=triangle*,mark repeat=3] table [x={Density}, y={Voronoi}] {Cook_TI_NC_2sinx_Vert.dat};
			\addlegendentry{Voronoi}
			\addplot[blue, solid,mark=diamond*,mark repeat=3] table [x={Density}, y={Q1}] {Cook_TI_NC_2sinx_Vert.dat};
			\addlegendentry{$Q_{1}$}
			\addplot[black, solid,mark=star,mark repeat=3] table [x={Density}, y={Q2}] {Cook_TI_NC_2sinx_Vert.dat};
			\addlegendentry{$Q_{2}$}
			\pgfplotsset{every axis/.append style={
					label style={font=\footnotesize},
					tick label style={font=\footnotesize}  
			}}
			\pgfplotsset{compat=1.3}
			\end{axis}
			\end{tikzpicture}
			\caption{}
		\end{subfigure}%
		\vskip\baselineskip
		\begin{subfigure}{0.5\textwidth}
			\centering
			\begin{tikzpicture}
			\begin{axis}[
			xlabel = {Mesh Density},
			ylabel = {Displacement},
			ymin = 1.4, ymax = 7.6,
			minor y tick num = 1,
			ytick = {1, 2, 3, 4, 5, 6, 7},
			every axis plot/.append style={thick},
			legend columns=2, 
			legend style={
				/tikz/column 2/.style={
					column sep=3.5pt}}]
			\pgfplotsset{every axis legend/.append style={	
					at={(0.97,0.04)},
					anchor=south east}}
			\addplot[magenta,mark=square*,mark repeat=3] table [x={Density}, y={VEM_quad}] {Cook_TI_NC_2sinx_w05.dat};
			\addlegendentry{Quad}
			\addplot[green,mark=pentagon*,mark repeat=3] table [x={Density}, y={VEM_hex}] {Cook_TI_NC_2sinx_w05.dat};
			\addlegendentry{Hex}
			\addplot[cyan,mark=triangle*,mark repeat=3] table [x={Density}, y={Voronoi}] {Cook_TI_NC_2sinx_w05.dat};
			\addlegendentry{Voronoi}
			\addplot[blue, solid,mark=diamond*,mark repeat=3] table [x={Density}, y={Q1}] {Cook_TI_NC_2sinx_w05.dat};
			\addlegendentry{$Q_{1}$}
			\addplot[black, solid,mark=star,mark repeat=3] table [x={Density}, y={Q2}] {Cook_TI_NC_2sinx_w05.dat};
			\addlegendentry{$Q_{2}$}
			\pgfplotsset{every axis/.append style={
					label style={font=\footnotesize},
					tick label style={font=\footnotesize}  
			}}
			\pgfplotsset{compat=1.3}
			\end{axis}
			\end{tikzpicture}
			\caption{}
		\end{subfigure}%
		\begin{subfigure}{0.5\textwidth}
			\centering
			\begin{tikzpicture}
			\begin{axis}[
			xlabel = {Mesh Density},
			ylabel = {Displacement},
			ymin = 1.4, ymax = 7.6,
			minor y tick num = 1,
			ytick = {1, 2, 3, 4, 5, 6, 7},
			every axis plot/.append style={thick},
			legend columns=2, 
			legend style={
				/tikz/column 2/.style={
					column sep=3.5pt}}]
			\pgfplotsset{every axis legend/.append style={	
					at={(0.97,0.04)},
					anchor=south east}}
			\addplot[magenta,mark=square*,mark repeat=3] table [x={Density}, y={VEM_quad}] {Cook_TI_NC_2sinx_VaryWeight.dat};
			\addlegendentry{Quad}
			\addplot[green,mark=pentagon*,mark repeat=3] table [x={Density}, y={VEM_hex}] {Cook_TI_NC_2sinx_VaryWeight.dat};
			\addlegendentry{Hex}
			\addplot[cyan,mark=triangle*,mark repeat=3] table [x={Density}, y={Voronoi}] {Cook_TI_NC_2sinx_VaryWeight.dat};
			\addlegendentry{Voronoi}
			\addplot[blue, solid,mark=diamond*,mark repeat=3] table [x={Density}, y={Q1}] {Cook_TI_NC_2sinx_VaryWeight.dat};
			\addlegendentry{$Q_{1}$}
			\addplot[black, solid,mark=star,mark repeat=3] table [x={Density}, y={Q2}] {Cook_TI_NC_2sinx_VaryWeight.dat};
			\addlegendentry{$Q_{2}$}
			\pgfplotsset{every axis/.append style={
					label style={font=\footnotesize},
					tick label style={font=\footnotesize}  
			}}
			\pgfplotsset{compat=1.3}
			\end{axis}
			\end{tikzpicture}
			\caption{}
		\end{subfigure}%
		\setlength{\belowcaptionskip}{-20pt}
		\caption{Tip displacement vs mesh density for Cook's membrane problem and with fibre directions defined by curves $2\sin x$, and using (a) the centroidal value of $\ba$; (b) the average of nodal values of $\ba$; (c) equal weighting of nodal and centroidal values; and (d) a varying weighted average as in \eqref{w_c} and \eqref{a_ave}}
		\label{Fig:Cook2Sinx}
	\end{figure}

	Figure \ref{Fig:Cook2Sinx} shows tip displacement as a function of mesh density for fibres corresponding to the sinusoidal distribution. In Figure \ref{Fig:Cook2Sinx} (a) we  calculate $ \boldsymbol{a}_{\rm ave}$ based on its value at the element centroid. Again we see good performance by the VEM elements, with an accuracy comparable to that of $Q_2$. In contrast to the results for polynomial variation in Figures \ref{Fig:CookPoly4} (a) and (b), the coarse-mesh behaviour is somewhat erratic, with a smooth dependence on mesh density only after $d \approx 25$. Figures \ref{Fig:Cook2Sinx} (c) and (d) present results for the cases in which, respectively, firstly an equal weighting of centroidal and nodal values of $\ba$ is used, and secondly, using \eqref{a_ave}, a varying weight is used. Similar behaviour is seen when compared with the results in Figure \ref{Fig:CookPoly4} (b), though the dependence on mesh density becomes smoother for coarser meshes at $d \approx 15$.
	\newpage 
	
	\bigskip
	
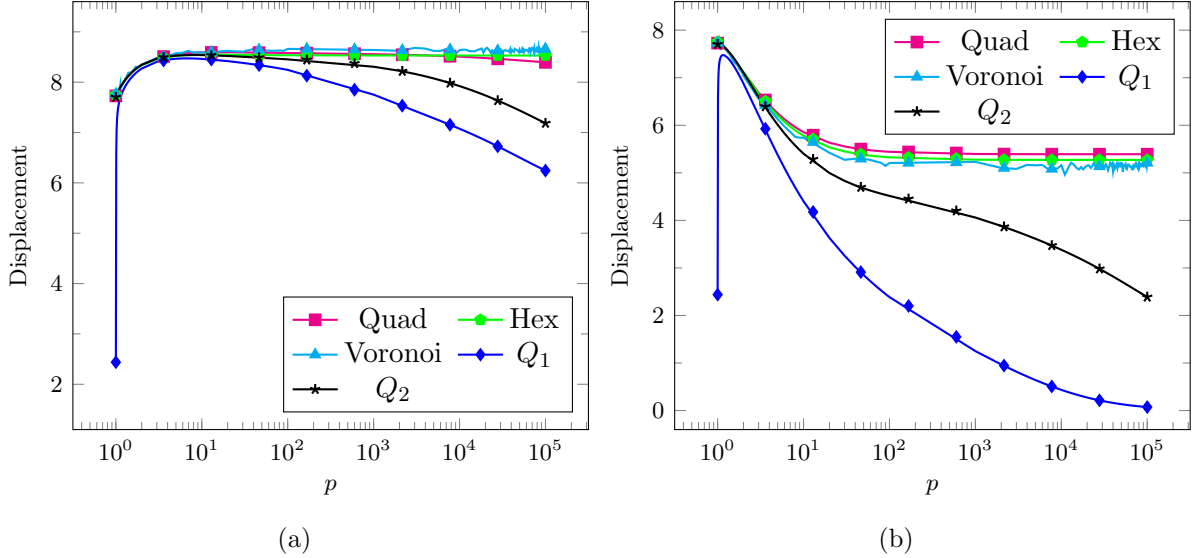
\begin{figure}[ht!]
		\centering
		\begin{subfigure}{0.5\textwidth}
			\centering
			\begin{tikzpicture}
			\begin{semilogxaxis}[
			xlabel = {$p$},
			ylabel = {Displacement},
			every axis plot/.append style={thick},
			ymin = 1.1, ymax = 9.6,
			minor y tick num = 1,
			legend columns=2, 
			legend style={
				/tikz/column 2/.style={
					column sep=3.5pt}}]
			\pgfplotsset{every axis legend/.append style={	
					at={(0.97,0.04)},
					anchor=south east}}
			\addplot[magenta,mark=square*] table [x={p}, y={VEM_quad}] {legend.dat};
			\addlegendentry{Quad};
			\addplot[green,mark=pentagon*] table [x={p}, y={VEM_hex}] {legend.dat};
			\addlegendentry{Hex};
			\addplot[cyan,mark=triangle*] table [x={p}, y={Voronoi}] {legend.dat};
			\addlegendentry{Voronoi};
			\addplot[blue,mark=diamond*] table [x={p}, y={Q1}] {legend.dat};
			\addlegendentry{$Q_{1}$};
			\addplot[black,mark=star] table [x={p}, y={Q2}] {legend.dat};
			\addlegendentry{$Q_{2}$};
			\addplot[magenta] table [x={p}, y={Quad}] {CookX4NIC.dat};
			\addplot[only marks,mark=square*,mark options={magenta}] table [x={p}, y={VEM_quad}] {CookX4NIC_m.dat};
			\addplot[green] table [x={p}, y={Hex}] {CookX4NIC.dat};
			\addplot[only marks,mark=pentagon*,mark options={green}] table [x={p}, y={VEM_hex}] {CookX4NIC_m.dat};
			\addplot[cyan] table [x={p}, y={Voronoi}] {CookX4NIC.dat};
			\addplot[only marks,mark=triangle*,mark options={cyan}] table [x={p}, y={Voronoi}] {CookX4NIC_m.dat};
			\addplot[blue, solid] table [x={p}, y={Q1}] {CookX4NIC.dat};
			\addplot[only marks,mark=diamond*,mark options={blue}] table [x={p}, y={Q1}] {CookX4NIC_m.dat};
			\addplot[black, solid] table [x={p}, y={Q2}] {CookX4NIC.dat};
			\addplot[only marks,mark=star,mark options={black}] table [x={p}, y={Q2}] {CookX4NIC_m.dat};
			%
			\pgfplotsset{every axis/.append style={
					label style={font=\footnotesize},
					tick label style={font=\footnotesize}  
			}}
			\pgfplotsset{compat=1.3}
			\end{semilogxaxis}
			\end{tikzpicture}
			\caption{ }
		\end{subfigure}%
		\begin{subfigure}{0.5\textwidth}
			\centering
			\begin{tikzpicture}
			\begin{semilogxaxis}[
			xlabel = {$p$},
			ylabel = {Displacement},
			every axis plot/.append style={thick},
			ymin = -0.4, ymax = 8.6,
			minor y tick num = 1,
			legend columns=2, 
			legend style={
				/tikz/column 2/.style={
					column sep=3.5pt}}]
			\pgfplotsset{every axis legend/.append style={	
					at={(0.97,0.96)},
					anchor=north east}}
			\addplot[magenta,mark=square*] table [x={p}, y={VEM_quad}] {legend.dat};
			\addlegendentry{Quad};
			\addplot[green,mark=pentagon*] table [x={p}, y={VEM_hex}] {legend.dat};
			\addlegendentry{Hex};
			\addplot[cyan,mark=triangle*] table [x={p}, y={Voronoi}] {legend.dat};
			\addlegendentry{Voronoi};
			\addplot[blue,mark=diamond*] table [x={p}, y={Q1}] {legend.dat};
			\addlegendentry{$Q_{1}$};
			\addplot[black,mark=star] table [x={p}, y={Q2}] {legend.dat};
			\addlegendentry{$Q_{2}$};
			\addplot[magenta] table [x={p}, y={Quad}] {Cook2sinxNIC.dat};
			\addplot[only marks,mark=square*,mark options={magenta}] table [x={p}, y={VEM_quad}] {Cook2sinxNIC_m.dat};
			\addplot[green] table [x={p}, y={Hex}] {Cook2sinxNIC.dat};
			\addplot[only marks,mark=pentagon*,mark options={green}] table [x={p}, y={VEM_hex}] {Cook2sinxNIC_m.dat};
			\addplot[cyan] table [x={p}, y={Voronoi}] {Cook2sinxNIC.dat};
			\addplot[only marks,mark=triangle*,mark options={cyan}] table [x={p}, y={Voronoi}] {Cook2sinxNIC_m.dat};
			\addplot[blue, solid] table [x={p}, y={Q1}] {Cook2sinxNIC.dat};
			\addplot[only marks,mark=diamond*,mark options={blue}] table [x={p}, y={Q1}] {Cook2sinxNIC_m.dat};
			\addplot[black, solid] table [x={p}, y={Q2}] {Cook2sinxNIC.dat};
			\addplot[only marks,mark=star,mark options={black}] table [x={p}, y={Q2}] {Cook2sinxNIC_m.dat};
			%
			\pgfplotsset{every axis/.append style={
					label style={font=\footnotesize},
					tick label style={font=\footnotesize}  
			}}
			\pgfplotsset{compat=1.3}
			\end{semilogxaxis}
			\end{tikzpicture}
			\caption{ }
		\end{subfigure}%
		\caption{Tip displacement vs $p$ for the Cook problem, for near-incompressibility and using (a) the quartic, and (b) the sinusoidal fibre distributions}
		\label{Cook_nic}
	\end{figure}

	Next, we consider behaviour in the near-incompressible limit, with $\nu = 0.49995$. Figures \ref{Cook_nic} (a) and (b) show tip displacement as a function of $p$, for the quartic and sinusoidal fibre distributions respectively. For the polynomial fibre distribution the centroidal values of fibre direction are used, while the weighted method is used for the sinusoidal distribution.
	There is little variation in the performance of the various VEM meshes, though for the Voronoi mesh and for near-inextensibility small scatter is observed. The sub-optimal behaviour of the $Q_2$ mesh seen in the Cook example in Figure \ref{Fig:Cook_DISPvsTHETA_Const} is not evident here. The $Q_1$ mesh again displays locking behaviour except in a narrow range of mild anisotropy. 
	
\begin{figure}[h]
		\begin{subfigure}{0.5\textwidth}
			\centering
			\includegraphics[width=0.8\textwidth]{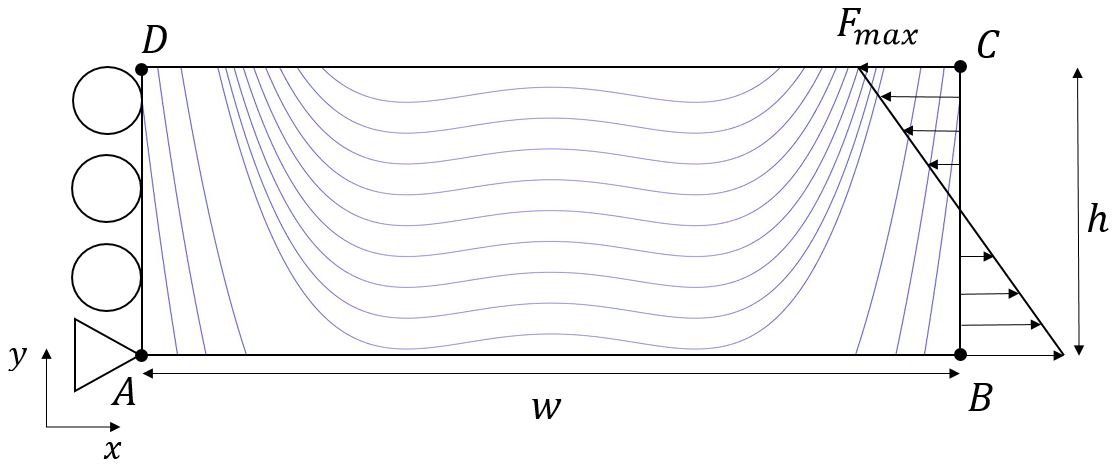}
		\end{subfigure}
		\begin{subfigure}{0.5\textwidth}
			\centering
			\includegraphics[width=0.8\textwidth]{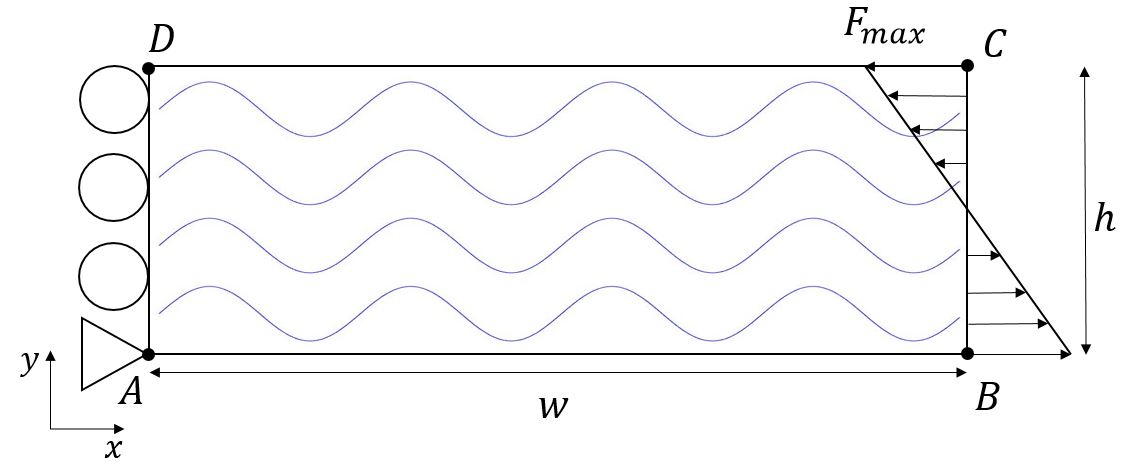}
		\end{subfigure}	
		\caption{Quartic and sinusoidal fibre distributions for the beam in bending problem}
		\label{Beam_quartic_2sin}
	\end{figure}

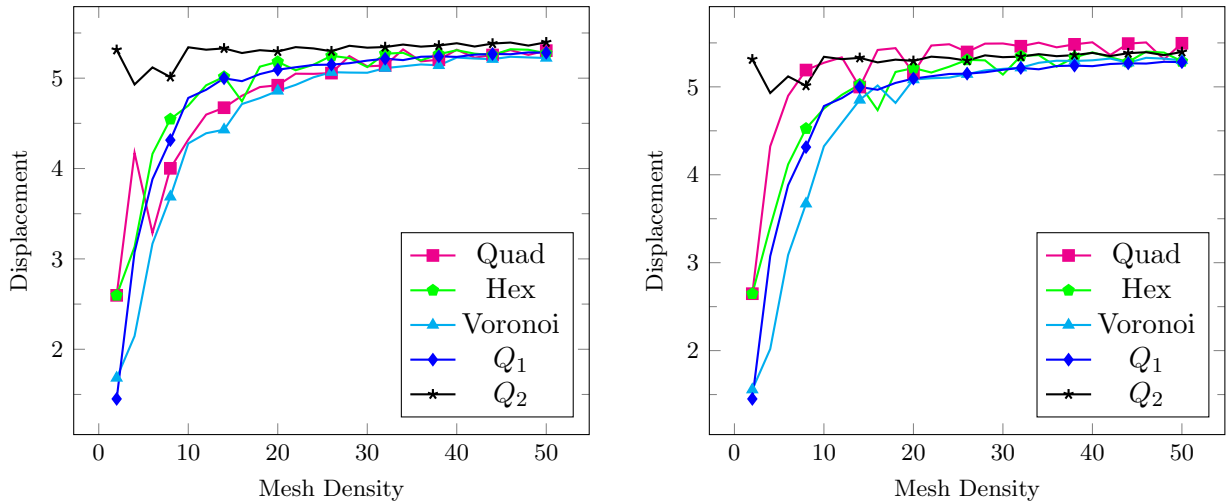
\begin{figure}[ht!]
		\centering
		\begin{subfigure}{0.45\textwidth}
			\centering
			\begin{tikzpicture}
			\begin{axis}[
			xlabel = {Mesh Density},
			ylabel = {Displacement},
			minor y tick num = 1,
			every axis plot/.append style={thick}]
			\addplot[magenta,mark=square*,mark repeat=3] table [x={Density}, y={VEM_quad}] {Beam_TI_NC_poly4_Vert.dat};
			\addlegendentry{Quad}
			\addplot[green,mark=pentagon*,mark repeat=3] table [x={Density}, y={VEM_hex}] {Beam_TI_NC_poly4_Vert.dat};
			\addlegendentry{Hex}
			\addplot[cyan,mark=triangle*,mark repeat=3] table [x={Density}, y={Voronoi}] {Beam_TI_NC_poly4_Vert.dat};
			\addlegendentry{Voronoi}
			\addplot[blue, solid,mark=diamond*,mark repeat=3] table [x={Density}, y={Q1}] {Beam_TI_NC_poly4_Vert.dat};
			\addlegendentry{$Q_{1}$}
			\addplot[black, solid,mark=star,mark repeat=3] table [x={Density}, y={Q2}] {Beam_TI_NC_poly4_Vert.dat};
			\addlegendentry{$Q_{2}$}
			\pgfplotsset{every axis legend/.append style={
					at={(0.97,0.04)},
					anchor=south east}}
			\pgfplotsset{every axis/.append style={
					label style={font=\footnotesize},
					tick label style={font=\footnotesize}  
			}}
			\pgfplotsset{compat=1.3}	
			\end{axis}
			\end{tikzpicture}
		\end{subfigure}
		\hspace{6ex}
		\begin{subfigure}{0.45\textwidth}
			\centering
			\begin{tikzpicture}
			\begin{axis}[
			xlabel = {Mesh Density},
			ylabel = {Displacement},
			minor y tick num = 1,
			every axis plot/.append style={thick}]
			\addplot[magenta,mark=square*,mark repeat=3] table [x={Density}, y={VEM_quad}] {Beam_TI_NC_poly4_Cent.dat};
			\addlegendentry{Quad}
			\addplot[green,mark=pentagon*,mark repeat=3] table [x={Density}, y={VEM_hex}] {Beam_TI_NC_poly4_Cent.dat};
			\addlegendentry{Hex}
			\addplot[cyan,mark=triangle*,mark repeat=3] table [x={Density}, y={Voronoi}] {Beam_TI_NC_poly4_Cent.dat};
			\addlegendentry{Voronoi}
			\addplot[blue, solid,mark=diamond*,mark repeat=3] table [x={Density}, y={Q1}] {Beam_TI_NC_poly4_Cent.dat};
			\addlegendentry{$Q_{1}$}
			\addplot[black, solid,mark=star,mark repeat=3] table [x={Density}, y={Q2}] {Beam_TI_NC_poly4_Cent.dat};
			\addlegendentry{$Q_{2}$}
			\pgfplotsset{every axis legend/.append style={
					at={(0.97,0.04)},
					anchor=south east}}
			\pgfplotsset{every axis/.append style={
					label style={font=\footnotesize},
					tick label style={font=\footnotesize}  
			}}
			\pgfplotsset{compat=1.3}	
			\end{axis}
			\end{tikzpicture}
		\end{subfigure}%
		\setlength{\belowcaptionskip}{-20pt}
		\caption{Tip displacement vs mesh density for the beam in bending problem, with $p=5$, and fibre directions defined by a quartic polynomial, using (a) the average of nodal values of $\ba$; (b) the centroidal value of $\ba$}
		\label{fig:BeamPoly4}
	\end{figure}

	\bigskip
	{\bf  Beam in bending.} We consider next the problem of a beam in bending, shown in Figure \ref{Beam_quartic_2sin}, with boundary conditions slightly different to those shown in Figure \ref{PureBendA_pi4_Fig}; the left edge is now constrained horizontally and pinned at the bottom left corner. The fibre distributions considered are once again quartic and sinusoidal, as for the Cook problem, and are depicted in Figure \ref{Beam_quartic_2sin}. However, the quartic distribution is now defined by $y = (x-5)^{2}(x-2.5)(x-7.5) + c$, and the sinusoidal distribution is as defined previously. The vertical displacement at point $ C $ is recorded. A value of Poisson's ratio of $\nu=0.3$ is used, except where indicted otherwise.
	\newline \noindent 
	In Figures \ref{fig:BeamPoly4} (a) and (b) we present convergence plots of vertical displacement vs mesh density, for mild anisotropy; that is, $p=5$, for the quartic distribution of fibres, and with $\ba_{\rm ave}$ calculated using respectively the centroidal value, and the average nodal value. We note stable convergence in all cases. For the case in which the average nodal direction is used, it is seen that the standard $Q_2$ element performs best for coarse meshes. 
	\newpage
	
\begin{figure}[ht!]
		\centering
		\begin{subfigure}{0.55\textwidth}
			\centering
			\begin{tikzpicture}
			\begin{axis}[
			xlabel = {Mesh Density},
			ylabel = {Displacement},
			minor y tick num = 1,
			every axis plot/.append style={thick},
			legend columns=2, 
			legend style={
				/tikz/column 2/.style={
					column sep=3.5pt}}]
			\pgfplotsset{every axis legend/.append style={	
					at={(0.97,0.04)},
					anchor=south east}}
			\addplot[magenta,mark=square*,mark repeat=3] table [x={Density}, y={VEM_quad}] {Beam_TI_NC_2sinx_Cent.dat};
			\addlegendentry{Quad}
			\addplot[green,mark=pentagon*,mark repeat=3] table [x={Density}, y={VEM_hex}] {Beam_TI_NC_2sinx_Cent.dat};
			\addlegendentry{Hex}
			\addplot[cyan,mark=triangle*,mark repeat=3] table [x={Density}, y={Voronoi}] {Beam_TI_NC_2sinx_Cent.dat};
			\addlegendentry{Voronoi}
			\addplot[blue, solid,mark=diamond*,mark repeat=3] table [x={Density}, y={Q1}] {Beam_TI_NC_2sinx_Cent.dat};
			\addlegendentry{$Q_{1}$}
			\addplot[black, solid,mark=star,mark repeat=3] table [x={Density}, y={Q2}] {Beam_TI_NC_2sinx_Cent.dat};
			\addlegendentry{$Q_{2}$}
			\pgfplotsset{every axis/.append style={
					label style={font=\footnotesize},
					tick label style={font=\footnotesize}  
			}}
			\pgfplotsset{compat=1.3}
			\end{axis}
			\end{tikzpicture}
			\caption{Fibre Direction At Centroid}
		\end{subfigure}%
		\begin{subfigure}{0.5\textwidth}
			\centering
			\begin{tikzpicture}
			\begin{axis}[
			xlabel = {Mesh Density},
			ylabel = {Displacement},
			minor y tick num = 1,
			every axis plot/.append style={thick},
			legend columns=2, 
			legend style={
				/tikz/column 2/.style={
					column sep=3.5pt}}]
			\pgfplotsset{every axis legend/.append style={	
					at={(0.97,0.04)},
					anchor=south east}}
			\addplot[magenta,mark=square*,mark repeat=3] table [x={Density}, y={VEM_quad}] {Beam_TI_NC_2sinx_Vert.dat};
			\addlegendentry{Quad}
			\addplot[green,mark=pentagon*,mark repeat=3] table [x={Density}, y={VEM_hex}] {Beam_TI_NC_2sinx_Vert.dat};
			\addlegendentry{Hex}
			\addplot[cyan,mark=triangle*,mark repeat=3] table [x={Density}, y={Voronoi}] {Beam_TI_NC_2sinx_Vert.dat};
			\addlegendentry{Voronoi}
			\addplot[blue, solid,mark=diamond*,mark repeat=3] table [x={Density}, y={Q1}] {Beam_TI_NC_2sinx_Vert.dat};
			\addlegendentry{$Q_{1}$}
			\addplot[black, solid,mark=star,mark repeat=3] table [x={Density}, y={Q2}] {Beam_TI_NC_2sinx_Vert.dat};
			\addlegendentry{$Q_{2}$}
			\pgfplotsset{every axis/.append style={
					label style={font=\footnotesize},
					tick label style={font=\footnotesize}  
			}}
			\pgfplotsset{compat=1.3}
			\end{axis}
			\end{tikzpicture}
			\caption{Average Fibre Direction At Vertices}
		\end{subfigure}%
		\vskip\baselineskip
		\begin{subfigure}{0.55\textwidth}
			\centering
			\begin{tikzpicture}
			\begin{axis}[
			xlabel = {Mesh Density},
			ylabel = {Displacement},
			minor y tick num = 1,
			every axis plot/.append style={thick},
			legend columns=2, 
			legend style={
				/tikz/column 2/.style={
					column sep=3.5pt}}]
			\pgfplotsset{every axis legend/.append style={	
					at={(0.97,0.04)},
					anchor=south east}}
			\addplot[magenta,mark=square*,mark repeat=3] table [x={Density}, y={VEM_quad}] {Beam_TI_NC_2sinx_w05.dat};
			\addlegendentry{Quad}
			\addplot[green,mark=pentagon*,mark repeat=3] table [x={Density}, y={VEM_hex}] {Beam_TI_NC_2sinx_w05.dat};
			\addlegendentry{Hex}
			\addplot[cyan,mark=triangle*,mark repeat=3] table [x={Density}, y={Voronoi}] {Beam_TI_NC_2sinx_w05.dat};
			\addlegendentry{Voronoi}
			\addplot[blue, solid,mark=diamond*,mark repeat=3] table [x={Density}, y={Q1}] {Beam_TI_NC_2sinx_w05.dat};
			\addlegendentry{$Q_{1}$}
			\addplot[black, solid,mark=star,mark repeat=3] table [x={Density}, y={Q2}] {Beam_TI_NC_2sinx_w05.dat};
			\addlegendentry{$Q_{2}$}
			\pgfplotsset{every axis/.append style={
					label style={font=\footnotesize},
					tick label style={font=\footnotesize}  
			}}
			\pgfplotsset{compat=1.3}
			\end{axis}
			\end{tikzpicture}
			\caption{Constant Weighted Average - $w=\frac{1}{2}$}
		\end{subfigure}%
		\begin{subfigure}{0.5\textwidth}
			\centering
			\begin{tikzpicture}
			\begin{axis}[
			xlabel = {Mesh Density},
			ylabel = {Displacement},
			minor y tick num = 1,
			every axis plot/.append style={thick},
			legend columns=2, 
			legend style={
				/tikz/column 2/.style={
					column sep=3.5pt}}]
			\pgfplotsset{every axis legend/.append style={	
					at={(0.97,0.04)},
					anchor=south east}}
			\addplot[magenta,mark=square*,mark repeat=3] table [x={Density}, y={VEM_quad}] {Beam_TI_NC_2sinx_VaryWeight.dat};
			\addlegendentry{Quad}
			\addplot[green,mark=pentagon*,mark repeat=3] table [x={Density}, y={VEM_hex}] {Beam_TI_NC_2sinx_VaryWeight.dat};
			\addlegendentry{Hex}
			\addplot[cyan,mark=triangle*,mark repeat=3] table [x={Density}, y={Voronoi}] {Beam_TI_NC_2sinx_VaryWeight.dat};
			\addlegendentry{Voronoi}
			\addplot[blue, solid,mark=diamond*,mark repeat=3] table [x={Density}, y={Q1}] {Beam_TI_NC_2sinx_VaryWeight.dat};
			\addlegendentry{$Q_{1}$}
			\addplot[black, solid,mark=star,mark repeat=3] table [x={Density}, y={Q2}] {Beam_TI_NC_2sinx_VaryWeight.dat};
			\addlegendentry{$Q_{2}$}
			\pgfplotsset{every axis/.append style={
					label style={font=\footnotesize},
					tick label style={font=\footnotesize}  
			}}
			\pgfplotsset{compat=1.3}
			\end{axis}
			\end{tikzpicture}
			\caption{Varying Weighted Average}
		\end{subfigure}
		\setlength{\belowcaptionskip}{-20pt}
		\caption{Tip displacement vs mesh density for the beam in bending problem, with fibre directions defined by curves $2\sin x$, and using (a) the centroidal value of $\ba$; (b) the average of nodal values of $\ba$; (c) equal weighting of nodal and centroidal values; and (d) a varying weighted average as in \eqref{w_c} and \eqref{a_ave}}
		\label{fig:Beam2Sinx}
	\end{figure}
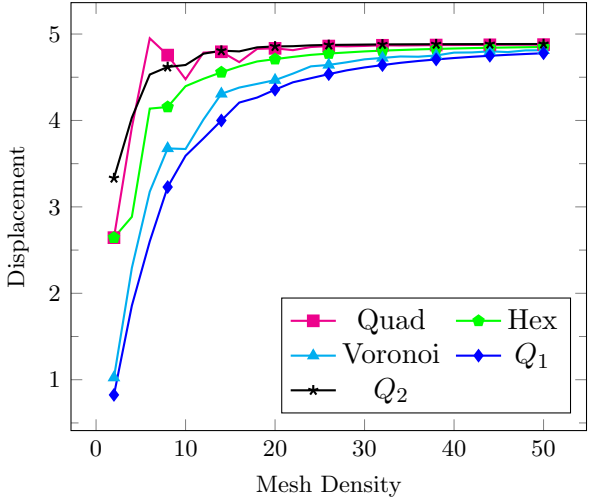
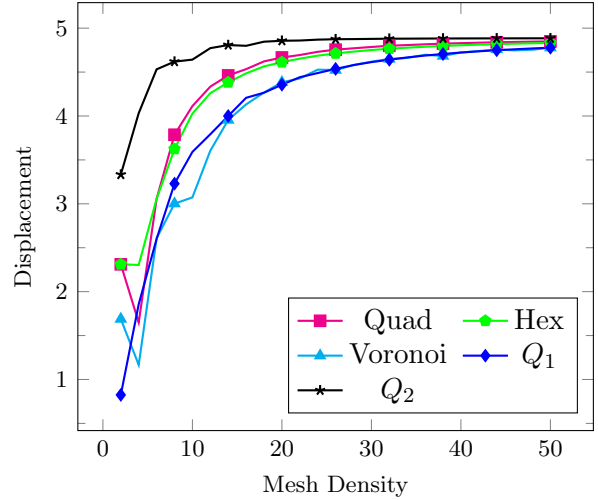
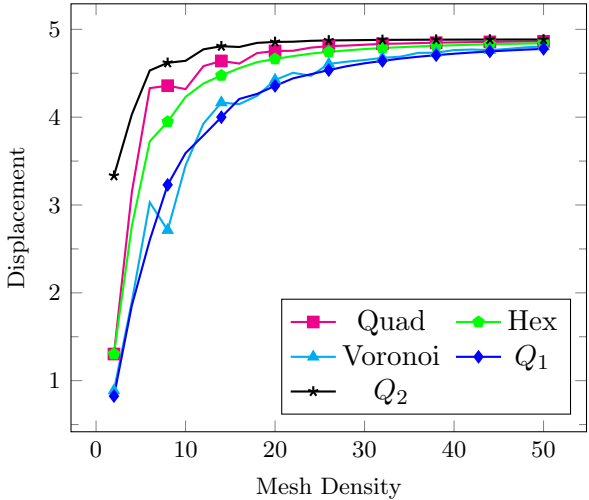
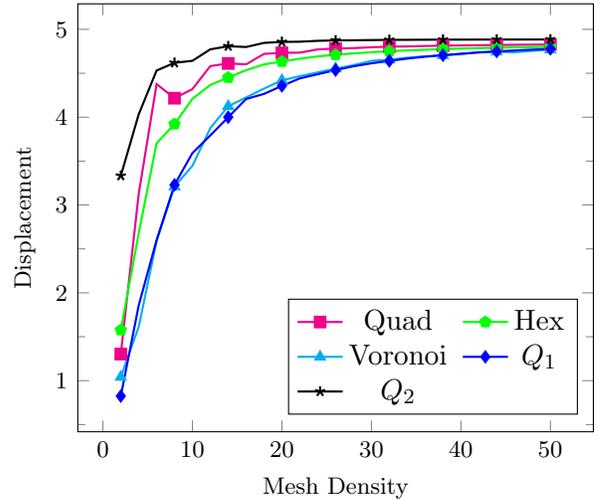

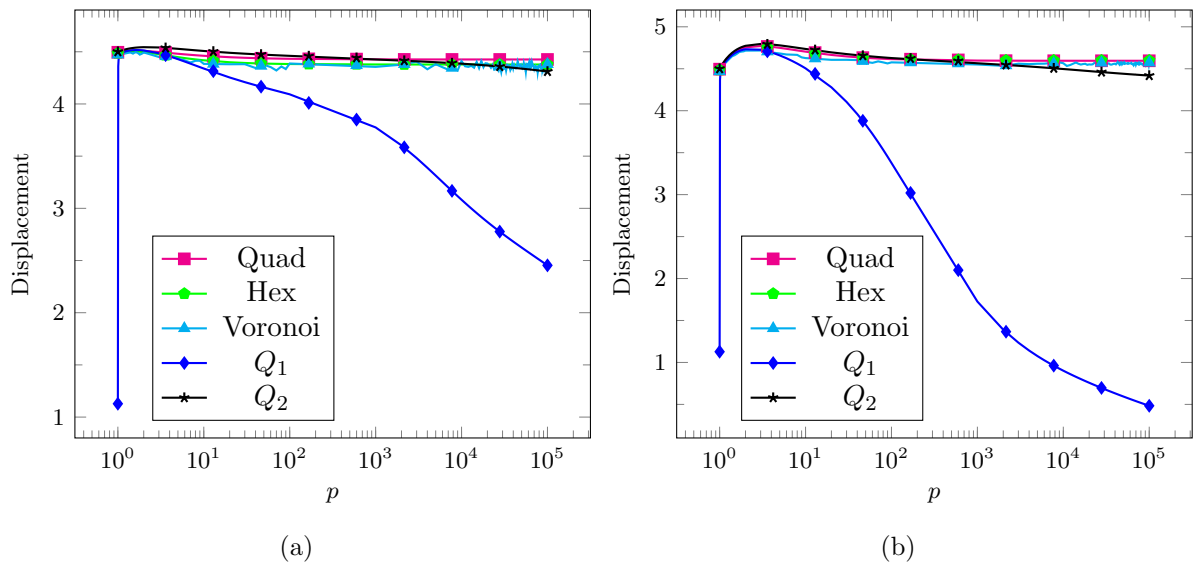
\begin{figure}[ht!]
		\centering
		\begin{subfigure}{0.5\textwidth}
			\centering
			\begin{tikzpicture}
			\begin{semilogxaxis}[
			xlabel = {$p$},
			ylabel = {Displacement},
			every axis plot/.append style={thick},
			ymin = 0.8, ymax = 4.9,
			minor y tick num = 1,
			legend columns=1, 
			legend style={
				/tikz/column 1/.style={
					column sep=3.5pt}}]
			\pgfplotsset{every axis legend/.append style={	
					at={(0.15,0.04)},
					anchor=south west}}
			\addplot[magenta,mark=square*] table [x={p}, y={VEM_quad}] {legend.dat};
			\addlegendentry{Quad};
			\addplot[green,mark=pentagon*] table [x={p}, y={VEM_hex}] {legend.dat};
			\addlegendentry{Hex};
			\addplot[cyan,mark=triangle*] table [x={p}, y={Voronoi}] {legend.dat};
			\addlegendentry{Voronoi};
			\addplot[blue,mark=diamond*] table [x={p}, y={Q1}] {legend.dat};
			\addlegendentry{$Q_{1}$};
			\addplot[black,mark=star] table [x={p}, y={Q2}] {legend.dat};
			\addlegendentry{$Q_{2}$};	
			\addplot[magenta] table [x={p}, y={Quad}] {BeamX4NIC.dat};
			\addplot[only marks,mark=square*,mark options={magenta}] table [x={p}, y={VEM_quad}] {BeamX4NIC_m.dat};
			\addplot[green] table [x={p}, y={Hex}] {BeamX4NIC.dat};
			\addplot[only marks,mark=pentagon*,mark options={green}] table [x={p}, y={VEM_hex}] {BeamX4NIC_m.dat};
			\addplot[cyan] table [x={p}, y={Voronoi}] {BeamX4NIC.dat};
			\addplot[only marks,mark=triangle*,mark options={cyan}] table [x={p}, y={Voronoi}] {BeamX4NIC_m.dat};
			\addplot[blue, solid] table [x={p}, y={Q1}] {BeamX4NIC.dat};
			\addplot[only marks,mark=diamond*,mark options={blue}] table [x={p}, y={Q1}] {BeamX4NIC_m.dat};
			\addplot[black, solid] table [x={p}, y={Q2}] {BeamX4NIC.dat};
			\addplot[only marks,mark=star,mark options={black}] table [x={p}, y={Q2}] {BeamX4NIC_m.dat};
			%
			\pgfplotsset{every axis/.append style={
					label style={font=\footnotesize},
					tick label style={font=\footnotesize}  
			}}
			\pgfplotsset{compat=1.3}
			\end{semilogxaxis}
			\end{tikzpicture}
			\caption{ }
		\end{subfigure}%
		\begin{subfigure}{0.5\textwidth}
			\centering
			\begin{tikzpicture}
			\begin{semilogxaxis}[
			xlabel = {$p$},
			ylabel = {Displacement},
			every axis plot/.append style={thick},
			ymin = 0.1, ymax = 5.2,
			minor y tick num = 1,
			legend columns=1, 
			legend style={
				/tikz/column 1/.style={
					column sep=3.5pt}}]
			\pgfplotsset{every axis legend/.append style={	
					at={(0.125,0.04)},
					anchor=south west}}
			\addplot[magenta,mark=square*] table [x={p}, y={VEM_quad}] {legend.dat};
			\addlegendentry{Quad};
			\addplot[green,mark=pentagon*] table [x={p}, y={VEM_hex}] {legend.dat};
			\addlegendentry{Hex};
			\addplot[cyan,mark=triangle*] table [x={p}, y={Voronoi}] {legend.dat};
			\addlegendentry{Voronoi};
			\addplot[blue,mark=diamond*] table [x={p}, y={Q1}] {legend.dat};
			\addlegendentry{$Q_{1}$};
			\addplot[black,mark=star] table [x={p}, y={Q2}] {legend.dat};
			\addlegendentry{$Q_{2}$};	
			\addplot[magenta] table [x={p}, y={Quad}] {Beam2sinxNIC.dat};
			\addplot[only marks,mark=square*,mark options={magenta}] table [x={p}, y={VEM_quad}] {Beam2sinxNIC_m.dat};
			\addplot[green] table [x={p}, y={Hex}] {Cook2sinxNIC.dat};
			\addplot[only marks,mark=pentagon*,mark options={green}] table [x={p}, y={VEM_hex}] {Beam2sinxNIC_m.dat};
			\addplot[cyan] table [x={p}, y={Voronoi}] {Beam2sinxNIC.dat};
			\addplot[only marks,mark=triangle*,mark options={cyan}] table [x={p}, y={Voronoi}] {Beam2sinxNIC_m.dat};
			\addplot[blue, solid] table [x={p}, y={Q1}] {Beam2sinxNIC.dat};
			\addplot[only marks,mark=diamond*,mark options={blue}] table [x={p}, y={Q1}] {Beam2sinxNIC_m.dat};
			\addplot[black, solid] table [x={p}, y={Q2}] {Beam2sinxNIC.dat};
			\addplot[only marks,mark=star,mark options={black}] table [x={p}, y={Q2}] {Beam2sinxNIC_m.dat};
			%
			\pgfplotsset{every axis/.append style={
					label style={font=\footnotesize},
					tick label style={font=\footnotesize}  
			}}
			\pgfplotsset{compat=1.3}
			\end{semilogxaxis}
			\end{tikzpicture}
			\caption{}
		\end{subfigure}%
		\caption{Tip displacement vs $p$ for the beam in bending problem, for near-incompressibility and using (a) the quartic, and (b) the sinusoidal fibre distributions}
		\label{beam_nic}
	\end{figure}

	\bigskip
	Figure \ref{fig:Beam2Sinx} shows tip displacement as a function of mesh density for a sinusoidal distribution of fibres. In Figure \ref{fig:Beam2Sinx} (a) we  calculate $ \boldsymbol{a}_{\rm ave}$ based on its centroidal value. Again we see good performance by the VEM elements, though the $Q_2$ mesh performs best. Rapid numerical convergence is however observed after a density of $d \approx 25$.
	\newpage
	\bigskip
	As with the Cook problem, we next consider behaviour in the near-incompressible limit with $\nu = 0.49995$. Figure \ref{beam_nic} shows tip displacement as a function of $p$, for the quartic and sinusoidal fibre distributions. For the polynomial fibre distribution the average value of fibre direction at the vertices is used, while the weighted method is used for the sinusoidal distribution. There is little variation in the performance of the various VEM meshes, though for the Voronoi mesh and for near-inextensibility small scatter is again observed. The $Q_2$ mesh performs rather poorly, displaying some evidence of mild locking. This should be compared with the sub-optimal behaviour seen in Figure \ref{Fig:Cook_DISPvsTHETA_Const} for constant fibre directions. Except in a narrow range of mild anisotropy, the $Q_1$ mesh displays locking behaviour. 
	
	%
	%
	
	\section{Concluding remarks}
	In this work we have formulated and implemented a Virtual Element Method for plane transversely isotropic elasticity, making provision for homogeneous as well as non-homogeneous bodies.  In the latter case, various options for taking account of the non-constant elasticity tensor are investigated. The formulations have been studied numerically through two model problems, and for three different kinds of polygonal meshes. The results have been compared against those obtained using conventional conforming finite element approximations with bilinear and biquadratic approximations.
	
	The VEM approximations are found to be locking-free for both near-incompressibility and 
	near-inextensibility, without the need to make modifications to the formulation. In the case of finite element approximations, the well-known volumetric locking behaviour of bilinear approximations is evident, except for a range of parameters corresponding to mild anisotropy. This behaviour is consistent with what has been shown in \cite{Rasolofoson-Grieshaber-Reddy2018}; for mild anisotropy the Lam\'e parameter related to the volumetric response is bounded. Locking does however occur in the inextensible limit. 
	
	There have been few studies of transverse isotropy in the context of development of new finite element and related methods. The present study and the work cited above constitute two new contributions. Further work is in progress on alternative formulations such as, for example, the use of discontinuous Galerkin methods. The extension to problems involving nonlinear material behaviour and large deformations is also in progress. It would be of interest to investigate the extension of the work presented here to include higher order VEMs as well as problems in three dimensions.
	
	\section*{Acknowledgement}
	
	This work was carried out with support from the National Research Foundation of South Africa, through the South African Research Chair in Computational Mechanics. The authors acknowledge with thanks this support.

\bibliographystyle{unsrt}
\bibliography{Lit_VEM_basics}

\begin{thebibliography}{10}

\bibitem{Belytschko-Liu-Moran2000}
Ted Belytschko, Wing~Kam Liu, Brian Moran, and Khalil Elkhodary.
\newblock {\em Nonlinear Finite Elements for Continua and Structures}.
\newblock Wiley, 2014.

\bibitem{Wriggers2008}
P.~Wriggers.
\newblock {\em {N}onlinear {F}inite {E}lement {Methods}}.
\newblock Springer, Berlin, Heidelberg, New York, 2008.

\bibitem{Boffi-Brezzi-Fortin2013}
D.~Boffi, F.~Brezzi, and M.~Fortin.
\newblock {\em Mixed {F}inite {E}lement {M}ethods and {A}pplications}.
\newblock Springer, New York, 2013.

\bibitem{Hughes1987}
T.J.R. Hughes.
\newblock {\em The {F}inite {E}lement {M}ethod. {L}inear {S}tatic and {D}ynamic
  {F}inite {E}lement {A}nalysis}.
\newblock Prentice-Hall, Englewood Cliffs, New Jersey, 1987.

\bibitem{Arnold-Brezzi-Cockburn-Marini2002}
D.N. Arnold, F.~Brezzi, B.~Cockburn, and L.D. Marini.
\newblock Unified analysis of discontinuous {G}alerkin methods for elliptic
  problems.
\newblock {\em SIAM Journal on Numerical Analysis}, 39:1749--1779, 2002.

\bibitem{Grieshaber-McBride-Reddy2015}
B.~Grieshaber, A.T. McBride, and B.D. Reddy.
\newblock Uniformly convergent interior penalty methods using multilinear
  approximations for problems in elasticity.
\newblock {\em SIAM Journal on Numerical Analysis}, 53:2255--2278, 2015.

\bibitem{Hansbo-Larson2002}
P.~Hansbo and M.G. Larson.
\newblock Discontinuous {G}alerkin methods for incompressible and nearly
  incompressible elasticity by {N}itsche's method.
\newblock {\em Computer Methods in Applied Mechanics and Engineering},
  191:1895--1908, 2002.

\bibitem{Wihler2004}
T.~P. Wihler.
\newblock Locking-free {DGFEM} for elasticity problems in polygons.
\newblock {\em {IMA} Journal of Numerical Analysis}, 24(1):45--75, jan 2004.

\bibitem{Beirao-etal2013a}
L.~{Beir\~{a}o da Veiga}, F.~Brezzi, A.~Cangiani, G.~Manzini, L.D. Marini, and
  A.~Russo.
\newblock Basic principles of virtual element methods.
\newblock {\em Mathematical Models and Methods in Applied Sciences},
  23(01):199--214, 2013.

\bibitem{Beirao-etal2014}
L.~{Beir\~{a}o da Veiga}, F.~Brezzi, L.D. Marini, and A.~Russo.
\newblock The hitchhiker{\rq}s guide to the virtual element method.
\newblock {\em Mathematical Models and Methods in Applied Sciences},
  24:1541--1573, 2014.

\bibitem{Gain-etal2014}
A.L. Gain, C.~Talischi, and G.H. Paulino.
\newblock On the virtual element method for three-dimensional linear elasticity
  problems on arbitrary polyhedral meshes.
\newblock {\em Computer Methods in Applied Mechanics and Engineering},
  282:132--160, 2014.

\bibitem{Chi-Beirao-Paulino2017}
H.~Chi, L.~Beir\~{a}o~da Veiga, and G.~Paulino.
\newblock Some basic formulations of the virtual element method ({VEM}) for
  finite deformations.
\newblock {\em Computer Methods in Applied Mechanics and Engineering},
  318:148--192, 2017.

\bibitem{Wriggers-Reddy-Rust-Hudobivnik2017}
P.~Wriggers, B.D. Reddy, W.~Rust, and B.~Hudobivnik.
\newblock Efficient virtual element formulations for compressible and
  incompressible finite deformations.
\newblock {\em Computational Mechanics}, 60:253--268, 2017.

\bibitem{Artioli-etal2017b}
E.~Artioli, L.~Beir\~{a}o~da Veiga, C.~Lovadina, and E.~Sacco.
\newblock Arbitrary order 2{D} virtual elements for polygonal meshes: part
  {II}, inelastic problem.
\newblock {\em Computational Mechanics}, 60:643--657, 2017.

\bibitem{Beirao-etal2015}
L.~Beir\~{a}o~da Veiga, C.~Lovadina, and D.~Mora.
\newblock A virtual element method for elastic and inelastic problems on
  polytope meshes.
\newblock {\em Computer Methods in Applied Mechanics and Engineering},
  295:327--346, 2015.

\bibitem{Wriggers-Hudobivnik2017}
P.~Wriggers and B.~Hudobivnik.
\newblock A low order virtual element formulation for finite elasto-plastic
  deformations.
\newblock {\em Computer Methods in Applied Mechanics and Engineering},
  327:459--477, 2017.

\bibitem{Wriggers-Rust-Reddy2016}
P.~Wriggers, W.~Rust, and B.D. Reddy.
\newblock A virtual element method for contact.
\newblock {\em Computational Mechanics}, 58:1039--1050, 2016.

\bibitem{Wriggers-Hudobivnik-Korelc2018}
P.~Wriggers, B.~Hudobivnik, and J.~Korelc.
\newblock Efficient low order virtual elements for anisotropic materials at
  finite strains.
\newblock In {\em Advances in Computational Plasticity}, Computational Methods
  in Applied Sciences, vol. 46. Springer, Berlin, 2018.

\bibitem{Auricchio-Scalet-Wriggers2017}
F.~Auricchio, G.~Scalet, and P.~Wriggers.
\newblock Fiber-reinforced materials: finite elements for the treatment of the
  inextensibility constraint.
\newblock {\em Computational Mechanics}, 60:905--922, 2017.

\bibitem{Rasolofoson-Grieshaber-Reddy2018}
Faraniaina Rasolofoson, BJ~Grieshaber, and B~Daya Reddy.
\newblock Finite element approximations for near-incompressible and
  near-inextensible transversely isotropic bodies.
\newblock {\em International Journal for Numerical Methods in Engineering,},
  Advanced online publication, 2018.

\bibitem{Exadaktylos2001}
G.E. Exadaktylos.
\newblock On the constraints and relations of elastic constants of transversely
  isotropic geomaterials.
\newblock {\em International Journal of Rock Mechanics and Mining Sciences},
  38:941--956, 2001.

\bibitem{Artioli-etal2017a}
E.~Artioli, L.~Beir\~{a}o~da Veiga, C.~Lovadina, and E.~Sacco.
\newblock Arbitrary order 2{D} virtual elements for polygonal meshes: part {I},
  elastic problem.
\newblock {\em Computational Mechanics}, 60:355--377, 2017.

\end{thebibliography}

\end{document}